# INTRODUCTORY TOPICS IN DISTRIBUTIONS OVER BINARY TEST FUNCTIONS

Serban E. Vlad
Str. Zimbrului, Nr. 3, Bl. PB68, Et.2, Ap.11, 3700,
Oradea, Romania E-mail: serbanvlad@excite.com

***Abstract*** *We note with* $\boldsymbol{B}_2$ *the Boole algebra with two elements. In the paper we take profit on the possibility of defining an integral-differential calculus for the* $\boldsymbol{R}^n \to \boldsymbol{B}_2$ *functions in order to construct the elements of a distributions theory over such test functions.*

## 1. Preliminaries

1.1 We note with $\boldsymbol{B}_2$ the set $\{0,1\}$ together with the discrete topology, the order $0 \le 1$ and the laws $'\oplus'$ the modulo 2 sum, respectively $'\cdot'$ the product.

1.2 We note with the same symbols $'\oplus'$, $'\cdot'$ respectively $'\cdot'$ the modulo 2 sum and the product of the binary functions, respectively the product of the binary functions with binary scalars.

1.3 The convergence of the sequence $w_n \in \boldsymbol{B}_2, n \in \boldsymbol{N}$ consists in the next property:

$$\exists a \in \boldsymbol{B}_2, \exists N \in \boldsymbol{N}, \forall n \ge N, w_n = a \tag{1}$$

If the (unique) limit $a$ is mentioned, then we have the usual notations $w_n \to a$ or

$$\lim_{n\to\infty} w_n = a \tag{2}$$

1.4 We note with $|\cdot|$ the number of elements of a finite set.

1.5 Let $w_i \in \boldsymbol{B}_2, i \in I$ a binary family, or a function $w : I \to \boldsymbol{B}_2$

$$w_i \overset{not}{=} w(i), i \in I \tag{1}$$

Its *support* is by definition the set:

$$supp\ w = \{i \mid w_i = 1\} = \{i \mid w(i) = 1\} \tag{2}$$

1.6 If the support of $w$ is finite, then we define the modulo 2 summation

$$\underset{i\in I}{\Xi} w_i = \begin{cases} 1, |\,supp\ w\,| \text{ is an odd number} \\ 0, |\,supp\ w\,| \text{ is an even number} \end{cases} \tag{1}$$

where $0$, the number of elements of $\emptyset$, is by definition an even number; if the support of $w$ is infinite, then $\underset{i\in I}{\Xi} w_i$ is the symbol of a generalized divergent series.

1.7 The *left* and the *right limits* $\varphi(t-0),\varphi(t+0)$ and the *left* and the *right derivatives* $D^-\varphi(t),D^+\varphi(t)$ of the function $\varphi:\mathbf{R}\to\mathbf{B}_2$, binary numbers ($t$ fixed), or binary functions ($t$ variable) are defined like this:

$$\exists t'<t,\forall\xi\in(t',t),\varphi(\xi)=\varphi(t-0),t\in\mathbf{R}\vee\{\infty\} \tag{1}$$

$$\exists t'>t,\forall\xi\in(t,t'),\varphi(\xi)=\varphi(t+0),t\in\{-\infty\}\vee\mathbf{R} \tag{2}$$

$$D^-\varphi(t)=\varphi(t-0)\oplus\varphi(t),t\in\mathbf{R} \tag{3}$$

$$D^+\varphi(t)=\varphi(t+0)\oplus\varphi(t),t\in\mathbf{R} \tag{4}$$

1.8 Other notations for the left limit and the right limit functions of $\varphi$ are $\varphi^-:\mathbf{R}\vee\{\infty\}\to\mathbf{B}_2$, respectively $\varphi^+:\{-\infty\}\vee\mathbf{R}\to\mathbf{B}_2$. We shall also note with $\varphi^-,\varphi^+$ the restrictions of the previous functions to $\mathbf{R}$.

1.9 The conditions of existence of the left limit and of left derivability, respectively of existence of the right limit and of right derivability coincide. If they are fulfilled under the form:

$$\varphi(t-0) \text{ and } \varphi(t+0) \text{ exist}$$

then $\varphi$ is called *differentiable in* $t\in\mathbf{R}$ and if $\varphi$ is differentiable in any $t$, then it is simply called *differentiable,* [3].

1.10 The space of the $\mathbf{R}\to\mathbf{B}_2$ differentiable functions is noted with $Diff$ or $Diff^{(1)}$ and it is a $\mathbf{B}_2$-algebra (relative to the sum of the functions '$\oplus$', the product of the functions '$\cdot$' and the product of the functions with scalars '$\cdot$').

1.11 We note with $\chi_A:\mathbf{R}\to\mathbf{B}_2$ the characteristic function of the set $A\subset\mathbf{R}$.

1.12 **Theorem** (*of representation of the* $\mathbf{R}\to\mathbf{B}_2$ *differentiable functions*) The next conditions are equivalent:

a) $\varphi\in Diff$

b) the families $t_z\in\mathbf{R},a_z,b_z\in\mathbf{B}_2,z\in\mathbf{Z}$ exist so that

b.1) $...<t_{-1}<t_0<t_1<...$

b.2) $\forall t',t''\in\mathbf{R},(t',t'')\wedge\{t_z\mid z\in\mathbf{Z}\}$ is finite

b.3) $$\varphi(t)=\underset{z\in\mathbf{Z}}{\Xi}(a_z\cdot\chi_{\{t_z\}}(t)\oplus b_{z+1}\cdot\chi_{(t_z,t_{z+1})}(t)) \tag{1}$$

1.13 In the previous theorem, a real family $\{t_z\mid z\in\mathbf{Z}\}$ that satisfies:

i) the condition b.1) is called *strictly increasing*

ii)the condition b.2) is called *locally finite*

1.14 Being given $\varphi\in Diff$ like at 1.12, the families $(t_z),(a_z),(b_z)$ are not unique, as we can easily see.

1.15 We show the way that the $Diff\to Diff$ operators of taking the limit and the derivative act on the functions $\varphi$ from 1.12 (1):

$$\varphi(t-0)=\underset{z\in\mathbf{Z}}{\Xi} b_{z+1}\cdot\chi_{(t_z,t_{z+1}]}(t) \tag{1}$$

$$D^{-}\varphi(t)=\underset{z\in\mathbf{Z}}{\Xi} b_{z+1}\cdot\chi_{(t_z,t_{z+1}]}(t)\oplus\underset{z\in\mathbf{Z}}{\Xi}(a_z\cdot\chi_{\{t_z\}}(t)\oplus b_{z+1}\cdot\chi_{(t_z,t_{z+1})}(t)) \quad (2)$$

$$=\underset{z\in\mathbf{Z}}{\Xi}(a_z\oplus b_z)\cdot\chi_{\{t_z\}}(t)=\underset{z\in\mathbf{Z}}{\Xi} D^{-}\varphi(t_z)\cdot\chi_{\{t_z\}}(t)$$

and similarly for $\varphi(t+0), D^{+}\varphi(t)$.

1.16 **Conclusions** a) There exist left-right dualities that will not always be treated in full details.

b) The operators of taking the limit and of derivation are linear

c)
$$\varphi((t-0)-0)=\varphi(t-0) \quad (1)$$
$$\varphi((t-0)+0)=\varphi(t+0) \quad (2)$$
$$D^{-}D^{-}=D^{-} \quad (3)$$
$$D^{-}D^{+}=D^{+} \quad (4)$$

d) The sets $supp\ D^{-}\varphi, supp\ D^{+}\varphi$ are locally finite.

1.17 We refer to 1.12 and we show the way that the $\mathbf{R}^2\to\mathbf{B}_2$ differentiable functions are represented (are defined). The next statements are equivalent by definition:

a) $\varphi_2\in Diff^{(2)}$ (the $\mathbf{B}_2$-algebra of the $\mathbf{R}^2\to\mathbf{B}_2$ differentiable functions)

b) the families $t_z,u_w\in\mathbf{R}, a_{zw},b_{zw},c_{zw},d_{zw}\in\mathbf{B}_2, z,w\in\mathbf{Z}$ exist so that

b.1) $\{t_z\mid z\in\mathbf{Z}\},\{u_w\mid w\in\mathbf{Z}\}$ are strictly increasing

b.2) $\{t_z\mid z\in\mathbf{Z}\},\{u_w\mid w\in\mathbf{Z}\}$ are locally finite

b.3) $$\varphi_2(t,u)=\underset{z\in\mathbf{Z}}{\Xi}\underset{w\in\mathbf{Z}}{\Xi}(a_{zw}\cdot\chi_{\{t_z\}}(t)\cdot\chi_{\{u_w\}}(u)\oplus b_{zw+1}\cdot\chi_{\{t_z\}}(t)\cdot\chi_{(u_w,u_{w+1})}(u)\oplus$$
$$\oplus c_{z+1w}\cdot\chi_{(t_z,t_{z+1})}(t)\cdot\chi_{\{u_w\}}(u)\oplus d_{z+1w+1}\cdot\chi_{(t_z,t_{z+1})}(t)\cdot\chi_{(u_w,u_{w+1})}(u)) \quad (1)$$

## 2. The Space of the Test Functions

2.1 Let $\varphi:\mathbf{R}\to\mathbf{B}_2$ and the next conditions:

a) $\varphi\in Diff$

b) $\varphi$ has a right limit in $-\infty$ and a left limit in $\infty$ and these limits are null:
$$\varphi(-\infty+0)=\varphi(\infty-0)=0 \quad (1)$$

c) there exist the real numbers $t'<t''$ so that $supp\ \varphi\subset[t',t'']$ i.e. $\varphi$ vanishes outside a compact set

d) $\varphi$ is of the form

$$\varphi(t)=\varphi(t_0)\cdot\chi_{\{t_0\}}(t)\oplus\varphi(\frac{t_0+t_1}{2})\cdot\chi_{(t_0,t_1)}(t)\oplus\varphi(t_1)\cdot\chi_{\{t_1\}}(t)\oplus\varphi(\frac{t_1+t_2}{2})\cdot\chi_{(t_1,t_2)}(t)\oplus... \quad (2)$$
$$...\oplus\varphi(\frac{t_k+t_{k+1}}{2})\cdot\chi_{(t_k,t_{k+1})}(t)\oplus\varphi(t_{k+1})\cdot\chi_{\{t_{k+1}\}}(t)$$

where $t\in\mathbf{R}, k\in\mathbf{N}$ and $t_0<t_1<...<t_{k+1}$

e) $supp\ \varphi$ belongs to the set ring $\mathbf{S}$ (relative to $\Delta,\wedge$) generated by the sets of the form $(t_0,t_1),\{t_2\}$ where $t_0,t_1,t_2\in\mathbf{R}$.

2.2 If $\varphi$ fulfills one of the next equivalent conditions from 2.1:

a) and b) $\Leftrightarrow$ a) and c) $\Leftrightarrow$ d) $\Leftrightarrow$ e)

then it is called *test function*, or *fundamental function.*

2.3 We note with $\mathbb{D}$ the set of the test functions (also called the fundamental space), organized from an algebraical point of view as $\boldsymbol{B}_2$-algebra.

2.4 The next property

$$\psi\in Diff,\varphi\in \mathbb{D}\Rightarrow \psi\cdot\varphi\in \mathbb{D} \quad (1)$$

shows that $\mathbb{D}$ is an ideal of $Diff$ .

2.5 We note for $\tau\in \boldsymbol{R}$ and $h:\boldsymbol{R}\to \boldsymbol{B}_2$ arbitrary:

$$h_\tau(t) \overset{def}{=} h(t-\tau),t\in \boldsymbol{R} \quad (1)$$

2.6 It is easily shown that $\varphi\in\mathbb{D}$ implies $\varphi_\tau,\varphi^-,\varphi^+\in\mathbb{D}$. Moreover, $supp\ \varphi\subset[t_0,t_{k+1}]$ implies $supp\ \varphi_\tau\subset[t_0+\tau,t_{k+1}+\tau]$ and $supp\ \varphi^-, supp\ \varphi^+\subset[t_0,t_{k+1}]$.

2.7 Let us consider the triple $(\psi,(\tau_n),\varphi)$, where $\psi\in Diff,(\tau_n)$ is a real positive sequence converging to $0$ and $\varphi\in\mathbb{D}$. It defines the following sequences of test functions: $(\psi\cdot\varphi_{\tau_n}),(\psi\cdot\varphi_{-\tau_n})$ that we shall call *left convergent*, respectively *right convergent*.

## 3. Regular Distributions

3.1 The $\boldsymbol{B}_2$-algebra $I_{Loc}$, or $I_{Loc}^{(1)}$ of the *locally integrable functions* consists in the functions $f:\boldsymbol{R}\to\boldsymbol{B}_2$ with a *locally finite support*:

$$\forall t',t''\in \boldsymbol{R},(t',t'')\wedge supp\ f \text{ is finite}$$

It is a subalgebra of $Diff$ .

3.2 The $\boldsymbol{B}_2$-algebra $I_\infty$ of the *integrable functions* consists in these functions $g:\boldsymbol{R}\to\boldsymbol{B}_2$ that have a finite support.

It is a subalgebra of $I_{Loc}$.

3.3 Let $g\in I_\infty$. We note by $\int_{-\infty}^{\infty} g$ the *integral* of $g$ defined in this way:

$$\int_{-\infty}^{\infty} g = \underset{\xi\in\boldsymbol{R}}{\Xi}\ g(\xi) \quad (1)$$

3.4 Let $f\in I_{Loc}$ a locally integrable function. It defines the function $[f]:\mathbb{D}\to\boldsymbol{B}_2$ in the following manner:

$$[f](\varphi)=\int_{-\infty}^{\infty} f\cdot\varphi,\ \varphi\in\mathbb{D} \quad (1)$$

3.5 **Remark** If $supp\ \varphi\subset[t_0,t_{k+1}]$ then

$$supp\ f\cdot\varphi = supp\ f\wedge supp\ \varphi\subset supp\ f\wedge[t_0,t_{k+1}] \quad (1)$$

is a finite set, thus the integral 3.4 (1) has sense.

3.6 Let us suppose that $\varphi$ is given by the formula 2.1 (2), where $t_0,t_1,...,t_{k+1}$ are chosen so that $supp\ f\wedge[t_0,t_{k+1}]\subset\{t_0,t_1,...,t_{k+1}\}$ (this is always possible without the loss of the generality). Then:

$$f(t)\cdot\varphi(t)=f(t_0)\cdot\varphi(t_0)\cdot\chi_{\{t_0\}}(t)\oplus f(t_1)\cdot\varphi(t_1)\cdot\chi_{\{t_1\}}(t)\oplus...$$
$$...\oplus f(t_{k+1})\cdot\varphi(t_{k+1})\cdot\chi_{\{t_{k+1}\}}(t) \quad (1)$$

where $t\in\boldsymbol{R}$ and it takes place

$$[f](\varphi)=f(t_0)\cdot\varphi(t_0)\oplus f(t_1)\cdot\varphi(t_1)\oplus...\oplus f(t_{k+1})\cdot\varphi(t_{k+1}) \quad (2)$$

3.7 **Theorem** The function $[f]$ satisfies:

a) it is linear

b) for any triple $(\psi,(\tau_n),\varphi)$, where $\psi\in Diff,(\tau_n)$ is a real positive sequence that converges to $0$ and $\varphi\in\mathbf{D}$, the binary sequences $([f](\psi\cdot\varphi_{\tau_n})),([f](\psi\cdot\varphi_{-\tau_n}))$ are convergent.

**Proof** b) We suppose that $\varphi$ is given by the formula 2.1 (2) and $supp\ f\wedge[t_0,t_{k+1}]\subset\{t_0,t_1,...,t_{k+1}\}$. Then $supp\ f\cdot\psi\wedge[t_0,t_{k+1}]\subset\{t_0,t_1,...,t_{k+1}\}$ also and there exists $N\in\boldsymbol{N}$ such that for any $n\geq N$ we have:

$$[f](\psi\cdot\varphi_{\tau_n})=f(t_1)\cdot\psi(t_1)\cdot\varphi(t_1-0)\oplus...\oplus f(t_{k+1})\cdot\psi(t_{k+1})\cdot\varphi(t_{k+1}-0) \quad (1)$$

$$[f](\psi\cdot\varphi_{-\tau_n})=f(t_0)\cdot\psi(t_0)\cdot\varphi(t_0+0)\oplus...\oplus f(t_k)\cdot\psi(t_k)\cdot\varphi(t_k+0) \quad (2)$$

3.8 A simple and important property of the function $[f]$ is the following: for any $t\in\boldsymbol{R},\chi_{\{t\}}$ is a test function and

$$[f](\chi_{\{t\}})=\underset{\xi\in\boldsymbol{R}}{\Xi} f(\xi)\cdot\chi_{\{t\}}(\xi)=f(t) \quad (1)$$

3.9 **Proposition** Let $f,f_1\in I_{Loc}$. If

$$\forall\varphi\in\mathbf{D},[f](\varphi)=[f_1](\varphi) \quad (1)$$

then $f$ and $f_1$ are equal.

**Proof** Let us suppose that there exists $t_0\in\boldsymbol{R}$ with

$$f(t_0)\neq f_1(t_0) \quad (2)$$

This fact implies, taking into account 3.8, that

$$[f](\chi_{\{t_0\}})=f(t_0)\neq f_1(t_0)=[f_1](\chi_{\{t_0\}}) \quad (3)$$

which contradicts (1).

3.10 We discuss the behavior of $[f]:\mathbf{D}\to\boldsymbol{B}_2$, where $f\in I_{Loc}$, relative to translations and let $\tau\in\boldsymbol{R}$ an arbitrary number. Because $f_\tau\in I_{Loc}$, it makes sense to speak about the function $[f_\tau]:\mathbf{D}\to\boldsymbol{B}_2$.

For any $\varphi\in\mathbf{D}$ the next equations are true:

$$[f_\tau](\varphi)=\underset{\xi\in\boldsymbol{R}}{\Xi} f_\tau(\xi)\cdot\varphi(\xi)=\underset{\xi\in\boldsymbol{R}}{\Xi} f(\xi-\tau)\cdot\varphi(\xi)=\underset{\xi'\in\boldsymbol{R}}{\Xi} f(\xi')\cdot\varphi(\xi'+\tau)= \quad (1)$$
$$=\underset{\xi'\in\boldsymbol{R}}{\Xi} f(\xi')\cdot\varphi_{-\tau}(\xi')=[f](\varphi_{-\tau})$$

3.11 We discuss the behavior of $[f]:\mathbf{D}\to\boldsymbol{B}_2$ relative to the product with a differentiable function. We see that for any differentiable function $\psi\in Diff$ and any test function $\varphi\in\mathbf{D}$:

- the function $\psi\cdot f$ is locally integrable (i.e. $I_{Loc}$ is an ideal of $Diff$ )
- the function $\psi\cdot\varphi$ is a test function (see 2.4)

so that we can define $[\psi\cdot f]:\mathbf{D}\to\boldsymbol{B}_2$. It is true:

$$[\psi\cdot f](\varphi)=\int_{-\infty}^{\infty}(\psi\cdot f)\cdot\varphi=\int_{-\infty}^{\infty}f\cdot(\psi\cdot\varphi)=[f](\psi\cdot\varphi) \qquad (1)$$

3.12 The value of the function $[f]$ in the 'point' $\varphi$ is noted by tradition with $<[f],\varphi>$, $<[f](t),\varphi(t)>$ (abusive notation !), or $([f],\varphi)$, respectively $([f](t),\varphi(t))$ instead of $[f](\varphi)$.

3.13 The functions $[f]:\mathbf{D}\to \boldsymbol{B}_2$, where $f\in I_{Loc}$ that have been discussed in this paragraph are called *regular distributions*, or *distributions of function type*. It is said that $[f]$ is the *regular distribution defined by the function* $f$, or that is *associated to the locally integrable function* $f$.

## 4. Distributions

4.1 It is called *distribution* a function $f:\mathbf{D}\to \boldsymbol{B}_2$ that fulfills the following conditions:

a) it is linear

b) for any triple $(\psi,(\tau_n),\varphi)$, where $\psi\in Diff,(\tau_n)$ is a real positive sequence that converges to $0$ and $\varphi\in\mathbf{D}$, the binary sequences $(f(\psi\cdot\varphi_{\tau_n})),(f(\psi\cdot\varphi_{-\tau_n}))$ are convergent.

4.2 By tradition, the value of a distribution $f$ in the 'point' $\varphi$ is noted with $<f,\varphi>$, $<f(t),\varphi(t)>$ (abusive notation, because $t\in\boldsymbol{R}$ cannot be the argument of $f$), or $(f,\varphi)$, $(f(t),\varphi(t))$ instead of $f(\varphi)$.

4.3 a) The definition 4.1 restates the properties of the regular distributions that were expressed in 3.7 and which were considered to be essential.

b) We have already a remarkable example of distributions, that is the regular distributions.

4.4 The definition 4.1 leaves open the possibility that the following assertion:

$$\exists g\in I_{Loc},\forall\varphi\in\mathbf{D},<f,\varphi>=<[g],\varphi>$$

should be not true, thus $f$ may be not a regular distribution. In this situation, the distribution $f$ is said to be *singular*.

4.3 The fact that some distributions may be identified with the locally integrable functions, existing however distributions that can be identified with no (locally integrable) function was the cause for which I.M. Gelfand and G.E. Shilow have called the distributions *generalized functions*. The two great mathematicians were working with real functions, but we shall show the existence of the generalized pseudo-boolean functions a little later.

4.6 In order to interpret the conditions of convergence 4.1 b) let us suppose for the beginning that in $(\psi,(\tau_n),\varphi)$,

$$\psi=1 \qquad (1)$$

(the constant function). In this situation we get for $f:\mathbf{D}\to \boldsymbol{B}_2$ the transposition of the condition of differentiability (see 1.9) of $x:\boldsymbol{R}\to \boldsymbol{B}_2$ written under the form:

$b_1$) for any couple $((\tau_n),t)$ in which $(\tau_n)$ is a real positive sequence that converges to $0$ and $t\in\boldsymbol{R}$ is a real number, the binary sequences $(x(t-\tau_n)),(x(t+\tau_n))$ are convergent.

Thus, when $\psi\in Diff$ is arbitrary, 4.1 b) expresses under a sufficiently general form the idea that the distributions are 'differentiable'.

4.7 The problem that we put is if, continuing the previous idea, the way that $x(t-\tau_n)\to x(t-0)$ we have also $<f,\psi\cdot\varphi_{\tau_n}>\to<f,\psi\cdot\varphi^->$ or $<f,\varphi_{\tau_n}>\to<f,\varphi^->$. We shall show that the answer is generally negative.

4.8 We define the *parity function* $\pi:\boldsymbol{N}\to\boldsymbol{B}_2$ by:

$$\pi(n)=\begin{cases}1, n\ is\ odd\\0, n\ is\ even\end{cases} \quad (1)$$

see again 1.6 for another definition of the same type.

4.9 It is easily shown the next property of $\pi$:

$$\pi(m+n)=\pi(m)\oplus\pi(n),\ m,n\in\boldsymbol{N} \quad (1)$$

4.10 **Proposition** The function $f:\mathbf{D}\to\boldsymbol{B}_2$ defined by:

$$<f,\chi_{(a_1,b_1)}\oplus...\oplus\chi_{(a_p,b_p)}\oplus\chi_{\{c_1,...,c_k\}}>=\pi(p+k) \quad (1)$$

is a distribution, where $p,k\in\boldsymbol{N}$ and where, by definition:

$$p=0\Rightarrow\chi_{(a_1,b_1)}\oplus...\oplus\chi_{(a_p,b_p)}=0$$

$$k=0\Rightarrow\chi_{\{c_1,...,c_k\}}=0$$

the null functions.

**Proof** $f$ is obviously linear (see also 4.9) and on the other hand the binary sequences $(<f,\psi\cdot\varphi_{\tau_n}>),(<f,\psi\cdot\varphi_{-\tau_n}>)$ become constant starting with a certain rank.

4.11 **Counterexample** for 4.7 We take

$$\varphi=\chi_{(0,1)} \quad (1)$$

$$\tau_n=\frac{1}{n+1} \quad (2)$$

$$\psi=1 \quad (3)$$

and we have for $f$ like in 4.10 (1)

$$\varphi_{\tau_n}=\chi_{(\frac{1}{n+1},1+\frac{1}{n+1})} \quad (4)$$

$$\varphi^-=\chi_{(0,1]} \quad (5)$$

$$<f,\chi_{(\frac{1}{n+1},1+\frac{1}{n+1})}>=\pi(1)=1\neq 0=\pi(2)=<f,\chi_{(0,1]}> \quad (6)$$

where $n\in\boldsymbol{N}$.

4.12 Let the distribution $f:\mathbf{D}\to\boldsymbol{B}_2$ and $\psi\in Diff$, $\varphi\in\mathbf{D}$. The condition 4.1 b) implies that the function $g:\boldsymbol{R}\to\boldsymbol{B}_2$ defined in the next manner:

$$g(t)=<f,\psi\cdot\varphi_t>,t\in\boldsymbol{R} \quad (1)$$

is differentiable in the origin. Since $\psi$, $\varphi$ are arbitrarily chosen, they can be replaced by their translations with $\tau$, whence it follows that $g$ is differentiable in any $\tau$.

4.13 **Theorem** Let the distribution $f:\mathbf{D}\to\boldsymbol{B}_2$. The next statements are true:

a) The function $g:\boldsymbol{R}\to\boldsymbol{B}_2$ defined by

$$g(t)=<f,\varphi_t>,t\in\boldsymbol{R} \quad (1)$$

is differentiable, where $\varphi\in\mathbf{D}$.

b) The functions $F_0,h:\boldsymbol{R}\to\boldsymbol{B}_2$,

$$F_0(t) = < f, \chi_{\{t\}} > \quad (2)$$

$$h(t) = < f, \chi_{(t_0 + t, t_1 + t)} >, t \in \boldsymbol{R} \quad (3)$$

are differentiable.

c) The functions $x, y : \boldsymbol{R} \to \boldsymbol{B}_2$,

$$x(t) = < f, \chi_{(t_0, t)} > \quad (4)$$

$$y(t) = < f, \chi_{(t, t_0)} >, t \in \boldsymbol{R} \quad (5)$$

are differentiable, $t_0 \in \boldsymbol{R}$.

**Proof** We show that the function $x$ given by (4) has a right limit in $t_1 \in \boldsymbol{R}$, where $t_1$ is arbitrary.

Case 1
$$t_1 < t_0 \quad (6)$$

For any real positive sequence $\tau_n \in (0, t_0 - t_1), n \in \boldsymbol{N}$ that converges to $0$, we have:

$$< f, \chi_{(t_0, t_1 + \tau_n)} > = < f, \chi_{\varnothing} > = < f, 0 > = 0, n \in \boldsymbol{N} \quad (7)$$

Case 2
$$t_1 \geq t_0 \quad (8)$$

We suppose against all reason that $x$ does not have a right limit in $t_1$, thus the next positive convergent to 0 sequences exist: $(\tau_n), (\nu_n)$ with

$$< f, \chi_{(t_0, t_1 + \tau_n)} > \to a \quad (9)$$

$$< f, \chi_{(t_0, t_1 + \nu_n)} > \to a \oplus 1 \quad (10)$$

where $a \in \boldsymbol{B}_2$. We take the numbers $t_1^{'}$ and $t_0^{'}$ so that:

$$t_1^{'} < t_0 \leq t_1 < t_0^{'} \quad (11)$$

and $N \in \boldsymbol{N}$ exists giving for $n > N$ that

$$(t_1^{'} + \tau_n, t_1 + \tau_n) \wedge (t_0, t_0^{'}) = (t_0, t_1 + \tau_n) \quad (12)$$

$$(t_1^{'} + \nu_n, t_1 + \nu_n) \wedge (t_0, t_0^{'}) = (t_0, t_1 + \nu_n) \quad (13)$$

from where, taking into account that 4.1 b) implies the independence of the limit on $(\tau_n), (\nu_n)$:

$$< f, \chi_{(t_0, t_0^{'})} \cdot \chi_{(t_1^{'} + \tau_n, t_1 + \tau_n)} > = < f, \chi_{(t_0, t_1 + \tau_n)} > \to b \quad (14)$$

$$< f, \chi_{(t_0, t_0^{'})} \cdot \chi_{(t_1^{'} + \nu_n, t_1 + \nu_n)} > = < f, \chi_{(t_0, t_1 + \nu_n)} > \to b \quad (15)$$

$b \in \boldsymbol{B}_2$. The relations (14), (15) compared with (9), (10) give a contradiction and we have that $x$ has a right limit in $t_1$.

It is shown similarly that $x$ has a left limit in $t_1$ and the differentiability of $y$ too.

4.14 Let the distributions $f, g : \mathbb{D} \to \boldsymbol{B}_2$. We define the function $f \oplus g : \mathbb{D} \to \boldsymbol{B}_2$ by:

$$\forall \varphi \in \mathbb{D}, < f \oplus g, \varphi > = < f, \varphi > \oplus < g, \varphi > \quad (1)$$

4.15 **Theorem** The function $f \oplus g$ is a distribution.

**Proof** The property of linearity of $f \oplus g$, as well as the convergence of $(< f \oplus g, \psi \cdot \varphi_{\tau_n} >), (< f \oplus g, \psi \cdot \varphi_{-\tau_n} >)$ when $\psi \in Diff$, $(\tau_n)$ is real positive convergent to 0 and $\varphi \in \mathbb{D}$ - are obvious.

4.16 **Remarks** a) The theorem 4.15 justifies the notation $< f \oplus g, \varphi >$ for the value of the function $f \oplus g$ in the 'point' $\varphi$.

b) In the special case of the regular distributions, for $f, g \in I_{Loc}$ it is true:

$$[f \oplus g] = [f] \oplus [g] \tag{1}$$

4.17 Let the distribution $f : \mathbf{D} \to \mathbf{B}_2$ and the function $\psi \in Diff$. The property 3.11 of the regular distributions suggests the next definition: the function $\psi \cdot f : \mathbf{D} \to \mathbf{B}_2$ is given by:

$$< \psi \cdot f, \varphi > = < f, \psi \cdot \varphi >, \varphi \in \mathbf{D} \tag{1}$$

The definition is correct, because for any $\varphi$, the function $\psi \cdot \varphi$ is a test function.

4.18 **Theorem** The function $\psi \cdot f$ is a distribution.
**Proof** The linearity, as well as the property of convergence are obvious.

4.19 **Remarks** a) The previous theorem justifies the notation $< \psi \cdot f, \varphi >$ for the value of the function $\psi \cdot f$ in the 'point' $\varphi$.

b) We have the special case at 4.17, when $\psi$ is the constant function, that is identified with a Boolean constant.

4.20 We note with $\mathbf{D}'$ the space of the $\mathbf{D} \to \mathbf{B}_2$ distributions that is from an algebraical point of view organized as $\mathbf{B}_2$-linear space and $Diff$ -module.

4.21 Let the distribution $f \in \mathbf{D}'$ and $\tau \in \mathbf{R}$. The property 3.10 of the regular distributions suggests the next definition: the function $f_\tau : \mathbf{D} \to \mathbf{B}_2$ is given by:

$$< f_\tau, \varphi > = < f, \varphi_{-\tau} >, \varphi \in \mathbf{D} \tag{1}$$

The definition is correct because for any $\varphi \in \mathbf{D}, \varphi_{-\tau}$ is a test function.

4.22 **Theorem** $f_\tau$ is a distribution.
**Proof** We make use of the fact that

$$(\varphi \oplus \varphi')_{-\tau} = \varphi_{-\tau} \oplus \varphi'_{-\tau} \tag{1}$$

where $\varphi, \varphi' \in \mathbf{D}, \tau \in \mathbf{R}$ etc.

4.23 **Remarks** a) The previous theorem justifies the notation $< f_\tau, \varphi >$ for the value of $f_\tau$ in $\varphi$.

b) We have by the replacement of $\tau$ with $-\tau$:

$$< f_{-\tau}, \varphi > = < f, \varphi_\tau >, \varphi \in \mathbf{D} \tag{1}$$

c) If $[f] \in \mathbf{D}'$ is a regular distribution, then the next property is fulfilled:

$$[f_\tau] = [f]_\tau \tag{2}$$

## 5. Examples of Distributions

5.1 The null function $0 : \mathbf{D} \to \mathbf{B}_2$ is a distribution, the null element of the $\mathbf{B}_2$-linear space, respectively of the $Diff$ - module $\mathbf{D}'$.

5.2 We define the function $\delta : \mathbf{R} \to \mathbf{B}_2$ in the next manner:

$$\delta(t) = \chi_{\{0\}}(t) = \begin{cases} 1, t = 0 \\ 0, t \neq 0 \end{cases} \tag{1}$$

5.3 The next locally integrable functions define regular distributions: $\delta, \delta_{t_0}$, $\delta_{t_0} \oplus \delta_{t_1} \oplus ... \oplus \delta_{t_n}$, $\underset{z \in \mathbf{Z}}{\Xi} \delta_{t_z} : \mathbf{R} \to \mathbf{B}_2$, where $\{t_z \mid z \in \mathbf{Z}\}$ is a locally finite set:

$$<[\delta],\varphi>=\varphi(0) \quad (1)$$

$$<[\delta_{t_0}],\varphi>=\varphi(t_0) \quad (2)$$

$$<[\delta_{t_0}\oplus\delta_{t_1}\oplus...\oplus\delta_{t_n}],\varphi>=\varphi(t_0)\oplus\varphi(t_1)\oplus...\oplus\varphi(t_n) \quad (3)$$

$$<[\underset{z\in \boldsymbol{Z}}{\Xi}\delta_{t_z}],\varphi>=\underset{z\in \boldsymbol{Z}}{\Xi}\varphi(t_z) \quad (4)$$

In (1),…,(4), $\varphi\in\mathbf{D}$.

5.4 The function $\delta^-_{t_0}:\mathbf{D}\to \boldsymbol{B}_2$,

$$<\delta^-_{t_0},\varphi>=\varphi(t_0-0) \quad (1)$$

where $t_0\in \boldsymbol{R}$ and $\varphi\in\mathbf{D}$ - defines a singular distribution.

**Proof** a) We show that $\delta^-_{t_0}$ is a distribution. The linearity being obvious, let us consider the function $\psi\in Diff$, the real positive sequence $(\tau_n)$ that converges to 0 and the function $\varphi\in\mathbf{D}$. We have:

$$<\delta^-_{t_0},\psi\cdot\varphi_{\tau_n}>=\psi(t_0-0)\cdot\varphi(t_0-\tau_n-0)\to\psi(t_0-0)\cdot\varphi(t_0-0) \quad (2)$$

$$<\delta^-_{t_0},\psi\cdot\varphi_{-\tau_n}>=\psi(t_0-0)\cdot\varphi(t_0+\tau_n-0)\to\psi(t_0-0)\cdot\varphi(t_0+0) \quad (3)$$

b) We show that $\delta^-_{t_0}$ is singular. Against all reason, let $f\in I_{Loc}$ so that

$$\forall\varphi\in\mathbf{D},<\delta^-_{t_0},\varphi>=<[f],\varphi> \quad (4)$$

Because of 4.13 c), the binary sequence $(<\delta^-_{t_0},\chi_{(t_0-\tau_n,t_0)}>)$ has a limit, i.e. it converges towards

$$\chi_{(t_0-\tau_n,t_0)}(t_0-0)=1 \quad (5)$$

On the other hand, we have:

$$\exists N,\forall n\geq N,<[f],\chi_{(t_0-\tau_n,t_0)}>=\int_{-\infty}^{\infty}f\cdot\chi_{(t_0-\tau_n,t_0)}=\int_{-\infty}^{\infty}0=0 \quad (6)$$

resulting a contradiction. $\delta^-_{t_0}$ is singular.

5.5 The next distributions are singular: $\delta^-,\delta^+,\delta^-_{t_0}\oplus\delta^-_{t_1}\oplus...\oplus\delta^-_{t_n},\delta^+_{t_0}\oplus\delta^+_{t_1}\oplus...\oplus\delta^+_{t_n}$, $\underset{z\in \boldsymbol{Z}}{\Xi}\delta^-_{t_z},\underset{z\in \boldsymbol{Z}}{\Xi}\delta^+_{t_z}\in\mathbf{D}'$, where $\{t_z\,|\,z\in\boldsymbol{Z}\}$ is locally finite. We write only the 'left' examples:

$$<\delta^-,\varphi>=\varphi(0-0) \quad (1)$$

$$<\delta^-_{t_0}\oplus\delta^-_{t_1}\oplus...\oplus\delta^-_{t_n},\varphi>=\varphi(t_0-0)\oplus\varphi(t_1-0)\oplus...\oplus\varphi(t_n-0) \quad (2)$$

$$<\underset{z\in \boldsymbol{Z}}{\Xi}\delta^-_{t_z},\varphi>\overset{def}{=}\underset{z\in \boldsymbol{Z}}{\Xi}<\delta^-_{t_z},\varphi>=\underset{z\in \boldsymbol{Z}}{\Xi}\varphi(t_z-0) \quad (3)$$

with $\varphi\in\mathbf{D}$.

5.6 The next distributions are singular:

$$<\delta^-_{t_0}\oplus\delta_{t_0},\varphi>=\varphi(t_0-0)\oplus\varphi(t_0)=D^-\varphi(t_0) \quad (1)$$

$$<\underset{z\in \boldsymbol{Z}}{\Xi}(\delta^-_{t_z}\oplus\delta_{t_z}),\varphi>\overset{def}{=}\underset{z\in \boldsymbol{Z}}{\Xi}<\delta^-_{t_z}\oplus\delta_{t_z},\varphi>=\underset{z\in \boldsymbol{Z}}{\Xi}D^-\varphi(t_z) \quad (2)$$

where $\{t_z \mid z \in \mathbf{Z}\}$ is locally finite and $\varphi \in \mathbf{D}$ - as well as their right duals.

5.7 The functions

$$< f, \varphi >= \int_{-\infty}^{\infty} D^{-}\varphi \quad (1)$$

$$< g, \varphi >= \int_{-\infty}^{\infty} D^{+}\varphi \quad (2)$$

define singular distributions, where $\varphi \in \mathbf{D}$.

**Proof** a) We show that $f$ from (1) is a distribution and we refer to the formula of representation 2.1 (2).

a.1) The integral $\int_{-\infty}^{\infty} D^{-}\varphi$ makes sense for any $\varphi$; this fact results from the inclusion $supp\ D^{-}\varphi \subset \{t_0, t_1, ..., t_{k+1}\}$.

a.2) The formula of definition (1) of $f$ shows that it is linear (both $\int_{-\infty}^{\infty}$ and $D^{-}$ are linear).

a.3) We show that 4.1 b) is true. Without loss, for $\psi \in Diff$ we shall suppose that $supp\ D^{-}\psi, supp\ D^{+}\psi \subset \{..., t_{-1}, t_0, t_1, ..., t_{k+1}, t_{k+2}, ...\}$, where $(t_z)$ is locally finite. We get for $(\tau_n)$ positive convergent to 0 the existence of some $N \in \mathbf{N}$ with

$$t_0 < t_0 + \tau_n < t_1 < t_1 + \tau_n < ... < t_{k+1} < t_{k+1} + \tau_n < t_{k+2}$$

$$supp\ D^{-}(\psi \cdot \varphi_{\tau_n}) \subset \{t_0 + \tau_n, t_1, t_1 + \tau_n, ..., t_{k+1}, t_{k+1} + \tau_n\}$$

true whenever $n \geq N$, showing the convergence of $(< f, \psi \cdot \varphi_{\tau_n} >)$. The convergence of $(< f, \psi \cdot \varphi_{-\tau_n} >)$ is similarly proved.

b) In order to prove the singularity of $f$, let us suppose that some $h \in I_{Loc}$ exists with

$$\forall \varphi \in \mathbf{D}, \int_{-\infty}^{\infty} D^{-}\varphi = \int_{-\infty}^{\infty} h \cdot \varphi \quad (4)$$

We take

$$\varphi = \chi_{(t_0, t_1]} \quad (5)$$

with the property that $supp\ h \wedge (t_0, t_1] = \emptyset$. It is true:

$$\int_{-\infty}^{\infty} D^{-}\varphi = \int_{-\infty}^{\infty} \delta_{t_1} = 1 \neq 0 = \int_{-\infty}^{\infty} 0 = \int_{-\infty}^{\infty} h \cdot \varphi \quad (6)$$

contradiction. $f$ is singular.

5.8 Let $\varphi \in \mathbf{D}$ like at 2.1 (2). The association

$$\mathbf{D} \ni \varphi \mapsto < f, \varphi >= \varphi(t_0) \oplus \varphi(\frac{t_0 + t_1}{2}) \oplus \varphi(t_1) \oplus ... \oplus \varphi(\frac{t_k + t_{k+1}}{2}) \oplus \varphi(t_{k+1}) \in \mathbf{B}_2 \quad (1)$$

defines a function $f : \mathbf{D} \to \mathbf{B}_2$ that is a singular distribution.

**Proof** The fact that the function is a distribution was already stated at 4.10. We show that $f$ is singular. We have:

$$\forall n \in \boldsymbol{N}, < f, \delta_{\frac{1}{n+1}} > = 1 \tag{2}$$

On the other hand, for any $g \in I_{Loc}$ we have the existence of an $N$ so that $n \geq N$ implies:

$$g \cdot \delta_{\frac{1}{n+1}} = 0 \tag{3}$$

$$< [g], \delta_{\frac{1}{n+1}} > = \int_{-\infty}^{\infty} g \cdot \delta_{\frac{1}{n+1}} = 0 \tag{4}$$

The contradiction that we have obtained shows that $f$ is singular.

## 6. The Fundamental Functions that Are Associated to a Distribution. The Support of a Distribution

6.1 Let the distribution $f \in \mathbb{D}'$. We define the function $F : \boldsymbol{R}^2 \to \boldsymbol{B}_2$ in the next manner:

$$F(t', t'') = < f, \chi_{(t', t'')} >, t', t'' \in \boldsymbol{R} \tag{1}$$

6.2 **Remarks** a) If $t' \geq t''$, then

$$F(t', t'') = < f, \chi_{\varnothing} > = < f, 0 > = 0 \tag{1}$$

b) The functions $F(t', \cdot), F(\cdot, t'') : \boldsymbol{R} \to \boldsymbol{B}_2$ are differentiable, see 4.13 c).

c) We recall that the function $F_0 : \boldsymbol{R} \to \boldsymbol{B}_2$ from 4.13 (2)

$$F_0(t) = < f, \chi_{\{t\}} >, t \in \boldsymbol{R} \tag{2}$$

is differentiable.

6.3 With the previous notations, each distribution $f$ is of the form:

$$< f, \varphi(t_0) \cdot \chi_{\{t_0\}} \oplus \varphi(\frac{t_0 + t_1}{2}) \cdot \chi_{(t_0, t_1)} \oplus ... \oplus \varphi(\frac{t_k + t_{k+1}}{2}) \cdot \chi_{(t_k, t_{k+1})} \oplus \varphi(t_{k+1}) \cdot \chi_{\{t_{k+1}\}} > = \tag{1}$$

$$= \varphi(t_0) \cdot F_0(t_0) \oplus \varphi(\frac{t_0 + t_1}{2}) \cdot F(t_0, t_1) \oplus ... \oplus \varphi(\frac{t_k + t_{k+1}}{2}) \cdot F(t_k, t_{k+1}) \oplus \varphi(t_{k+1}) \cdot F_0(t_{k+1})$$

where $\varphi \in \mathbb{D}$.

6.4 The functions $F_*, F^* : \boldsymbol{R} \to \boldsymbol{B}_2$ are defined like this:

$$F_*(t) = \lim_{\substack{\varepsilon \to 0 \\ \varepsilon > 0}} < f, \chi_{(t - \varepsilon, t)} > \tag{1}$$

$$F^*(t) = \lim_{\substack{\varepsilon \to 0 \\ \varepsilon > 0}} < f, \chi_{(t, t + \varepsilon)} >, t \in \boldsymbol{R} \tag{2}$$

6.5 **Theorem** $F_*, F^*$ are differentiable.

**Proof** Let $t'' \in \boldsymbol{R}$ arbitrary and fixed and let $t' < t''$ chosen so that for any $\xi \in (t', t'')$:

$$F(\xi, t'') = a \tag{1}$$

$$F_0(\xi) = b \tag{2}$$

with $a, b \in \boldsymbol{B}_2$. Such $t', a, b$ exist because of 6.2 b) and c). If $\xi, \xi'$ satisfy $t' < \xi < \xi' < t''$, we have:

$$< f, \chi_{(\xi, t'')} > = < f, \chi_{(\xi, \xi')} > \oplus < f, \chi_{\{\xi'\}} > \oplus < f, \chi_{(\xi', t'')} > \tag{3}$$

$$a = < f, \chi_{(\xi, \xi')} > \oplus b \oplus a \tag{4}$$

and this is equivalent to the existence of $F_*(t''-0)$, i.e.

$$F_*(t''-0)=b=F_0(t''-0) \qquad (5)$$

As at the right of $t''$, $F(\cdot,t'')$ is null, it results that $F_*\in Diff$, because $t''$ was arbitrarily chosen.

In a similar manner it is shown that $F^*\in Diff$.

6.6 The four functions $F,F_0,F_*,F^*$ that were previously defined are called the *fundamental functions* that are associated to the distribution $f$.

6.7 The fundamental functions are not independent on each other. We have properties of the type:

a) $$F_0(t)=F(t',t)\oplus F(t,t'')\oplus F(t',t''), t'<t<t'' \qquad (1)$$

b) $$F_*(t)=\lim_{\substack{\varepsilon\to 0\\ \varepsilon>0}} F(t-\varepsilon,t) \qquad (2)$$

$$F^*(t)=\lim_{\substack{\varepsilon\to 0\\ \varepsilon>0}} F(t,t+\varepsilon) \qquad (3)$$

c) there exists the strictly increasing locally finite family $\{t_z \mid z\in \boldsymbol{Z}\}$ with

$$\forall z\in \boldsymbol{Z}, \forall t\in(t_z,t_{z+1}), F(t_z,t_{z+1})=F^*(t_z)\oplus F_0(t)\oplus F_*(t_{z+1}) \qquad (4)$$

**Proof** We prove c). Let $t_0\in \boldsymbol{R}$ be fixed; $t_1>t_0$ exists with the property that

$$F(t_0,t)=F^*(t_0) \qquad (5)$$

$$F_0(t)=F_0(t_0+0), t\in(t_0,t_1) \qquad (6)$$

from where the fact that

$$F(t,t_1)=F(t_0,t_1)\oplus F_0(t)\oplus F(t_0,t)\ (=F(t_0,t_1)\oplus F_0(t_0+0)\oplus F^*(t_0)) \qquad (7)$$

gives the conclusion

$$F_*(t_1)=F(t_0,t_1)\oplus F_0(t)\oplus F^*(t_0) \qquad (8)$$

i.e. $t_0,t_1$ are like at (4).

Let us suppose that $t_0<t_1<...<t_z$ are already defined. We take $t_{z+1}>t_z$ in the following way:

$$F(t_z,t)=F^*(t_z) \qquad (9)$$

$$F_0(t)=F_0(t_z+0), t\in(t_z,t_{z+1}) \qquad (10)$$

resulting

$$F(t,t_{z+1})=F(t_z,t_{z+1})\oplus F_0(t)\oplus F(t_z,t)\ (=F(t_z,t_{z+1})\oplus F_0(t_z+0)\oplus F^*(t_z)) \qquad (11)$$

and

$$F_*(t_{z+1})=F(t_z,t_{z+1})\oplus F_0(t)\oplus F^*(t_z) \qquad (12)$$

i.e. $t_0,...,t_{z+1}$ act like in (4).

The sequence $(t_z)_{z\geq 0}$ that is obtained iteratively may be chosen to be strictly increasing locally finite, because the supposition that a convergent subsequence $(t_{z_p})_p$ necessarily exists gives a contradiction. In a similar manner, there are obtained iteratively $t_{-1},t_{-2},...$

6.8 a) The distribution $f \in \mathbf{D}'$ defines a function $F : \boldsymbol{R}^2 \to \boldsymbol{B}_2$ by the relation 6.1 (1). $F$ satisfies the properties 6.2:

i) $t' \geq t'' \Rightarrow F(t', t'') = 0$ (1)

ii) $F(t', \cdot), F(\cdot, t'') \in Diff, t', t'' \in \boldsymbol{R}$

b) Conversely, let us remark for the beginning that any $\varphi \in \mathbf{D}$ may be written as a sum of functions of the form $\chi_{(t',t'')}$, if we take into account that

$$\chi_{\{t\}} = \chi_{(t',t)} \oplus \chi_{(t,t'')} \oplus \chi_{(t',t'')}, t' < t < t'' \quad (2)$$

Let now $F : \boldsymbol{R}^2 \to \boldsymbol{B}_2$ be given so that i), ii) are fulfilled. We define the function $f : \mathbf{D} \to \boldsymbol{B}_2$ by:

$$< f, \chi_{(t_1', t_1'')} \oplus ... \oplus \chi_{(t_k', t_k'')} > = F(t_1', t_1'') \oplus ... \oplus F(t_k', t_k''), k \geq 1 \quad (3)$$

$t_1', t_1'', ..., t_k', t_k'' \in \boldsymbol{R}$ and this represents the extension of the formula 6.1 (1) at all the test functions $\varphi \in \mathbf{D}$, by linearity. We show that $f \in \mathbf{D}'$, the linearity being obvious.

Let $(\psi, (\tau_n), \varphi)$ so that $\psi \in Diff$ satisfies (without loss)

$supp\ D^- \psi, supp\ D^+ \psi \subset \{..., t_{-1}, t_0, t_1, ..., t_{k+1}, ...\}$, $(\tau_n)$ is a real positive sequence convergent to 0 and $\varphi \in \mathbf{D}$ is given by 2.1 (2). There exists $N \in \boldsymbol{N}$ so that $n \geq N$ implies:

$$... < t_{-1} < t_0 < t_0 + \tau_n < t_1 < t_1 + \tau_n < ... < t_{k+1} < t_{k+1} + \tau_n < t_{k+2} < ...$$

$$\psi \cdot \varphi_{\tau_n} = \psi(\frac{t_0 + t_1}{2}) \cdot \varphi(t_0) \cdot \chi_{\{t_0 + \tau_n\}} \oplus \psi(\frac{t_0 + t_1}{2}) \cdot \varphi(\frac{t_0 + t_1}{2}) \cdot \chi_{(t_0 + \tau_n, t_1)} \oplus \quad (4)$$

$$\oplus \psi(t_1) \cdot \varphi(\frac{t_0 + t_1}{2}) \cdot \chi_{\{t_1\}} \oplus \psi(\frac{t_1 + t_2}{2}) \cdot \varphi(\frac{t_0 + t_1}{2}) \cdot \chi_{(t_1, t_1 + \tau_n)} \oplus$$

$$\oplus \psi(\frac{t_1 + t_2}{2}) \cdot \varphi(t_1) \cdot \chi_{\{t_1 + \tau_n\}} \oplus ... \oplus \psi(\frac{t_{k+1} + t_{k+2}}{2}) \cdot \varphi(t_{k+1}) \cdot \chi_{\{t_{k+1} + \tau_n\}}$$

the convergence of $(< f, \psi \cdot \varphi_{\tau_n} >)$ resulting from the linearity, from ii) and from (2).

In a similar manner the convergence of $(< f, \psi \cdot \varphi_{-\tau_n} >)$ is proved. $f$ is a distribution.

6.9 The *support* of $f \in \mathbf{D}'$ is by definition the support of the function $F : \boldsymbol{R}^2 \to \boldsymbol{B}_2$,

$$supp\ f \overset{def}{=} supp\ F = \{(t', t'') \mid (t', t'') \in \boldsymbol{R}^2, F(t', t'') = 1\} \quad (1)$$

6.10 We have a one-one association between the distributions $f \in \mathbf{D}'$ and the functions $F : \boldsymbol{R}^2 \to \boldsymbol{B}_2$ having the properties 6.8 i), ii). In a similar manner to the one-one association between the functions $x : \boldsymbol{R} \to \boldsymbol{B}_2$ and their supports $supp\ x \subset \boldsymbol{R}$, we have also a one-one association between the distributions $f \in \mathbf{D}'$ and their supports.

6.11 **Theorem** The distribution $f \in \mathbf{D}'$ is regular if and only if the next conditions are true:

i) $F_0 \in I_{Loc}$

ii) $F_*(t) = F^*(t) = 0, t \in \boldsymbol{R}$ (1)

**Proof** Only if Let us suppose that $[f]$ is a regular distribution, $f \in I_{Loc}$. Then

$$F_0(t) = < [f], \chi_{\{t\}} > = \int_{-\infty}^{\infty} f \cdot \chi_{\{t\}} = f(t) \in I_{Loc}, t \in \boldsymbol{R} \quad (2)$$

and on the other hand, for any $t\in \boldsymbol{R}$, there exists $t'<t$ so that

$$supp\ f\wedge(t',t)=\varnothing \tag{3}$$

giving

$$F_*(t)=\int_{-\infty}^{\infty} f\cdot\chi_{(t',t)}=\int_{-\infty}^{\infty}0=0 \tag{4}$$

and similarly for $F^*(t)$.

If Because $F_0\in I_{Loc}$, it makes sense to refer to the distribution $[F_0]\in \mathbf{D}'$ and let $\varphi\in\mathbf{D}$.

There exists the strictly increasing locally finite family $\{t_z \mid z\in \boldsymbol{Z}\}$ (see 6.7 c)) with the property that $\varphi$ can be put under the form:

$$\varphi(t)=\chi_{(a_1,b_1)}(t)\oplus...\oplus\chi_{(a_p,b_p)}(t)\oplus\chi_{\{c_1,...,c_n\}}(t),\ t\in \boldsymbol{R} \tag{5}$$

$p,n\in \boldsymbol{N}$ and the following conditions are fulfilled:

i) $a_1,b_1,...,a_p,b_p,c_1,...,c_n\in\{t_z \mid z\in \boldsymbol{Z}\}$ (if $p=0$ or $n=0$, then the condition is superfluous for these terms)

ii) $$\forall z\in \boldsymbol{Z},\forall t\in(t_z,t_{z+1}),F(t_z,t_{z+1})=F_0(t)=0 \tag{6}$$

It results that

$$\begin{aligned}<f,\varphi>&=<f,\chi_{(a_1,b_1)}\oplus...\oplus\chi_{(a_p,b_p)}\oplus\chi_{\{c_1,...,c_n\}}>= \\ &=<f,\chi_{(a_1,b_1)}>\oplus...\oplus<f,\chi_{(a_p,b_p)}>\oplus<f,\chi_{\{c_1\}}>\oplus...\oplus<f,\chi_{\{c_n\}}>= \\ &=F(a_1,b_1)\oplus...\oplus F(a_n,b_n)\oplus F_0(c_1)\oplus...\oplus F_0(c_n)= \\ &=F_0(c_1)\oplus...\oplus F_0(c_n)=\int_{-\infty}^{\infty}F_0\cdot\varphi=<[F_0],\varphi>\end{aligned} \tag{7}$$

## 7. The Limits and the Derivatives of the Distributions

7.1 We recall that for any real positive sequence $(\tau_n)$ convergent to 0 and for any test function $\varphi\in\mathbf{D}$, the distribution $f\in\mathbf{D}'$ being arbitrary, the binary sequences $(<f,\varphi_{\tau_n}>)$, $(<f,\varphi_{-\tau_n}>)$ are convergent and this means the existence of the next limits:

$$\lim_{\substack{\varepsilon\to 0\\ \varepsilon>0}}<f_{-\varepsilon},\varphi>=\lim_{\substack{\varepsilon\to 0\\ \varepsilon>0}}<f,\varphi_\varepsilon> \tag{1}$$

$$\lim_{\substack{\varepsilon\to 0\\ \varepsilon>0}}<f_{\varepsilon},\varphi>=\lim_{\substack{\varepsilon\to 0\\ \varepsilon>0}}<f,\varphi_{-\varepsilon}> \tag{2}$$

7.2 **Notation** Let $f\in\mathbf{D}'$. The functions $f^-,f^+:\mathbf{D}\to \boldsymbol{B}_2$ are defined in the next manner:

$$<f^-,\varphi>=\lim_{\substack{\varepsilon\to 0\\ \varepsilon>0}}<f,\varphi_\varepsilon> \tag{1}$$

$$<f^+,\varphi>=\lim_{\substack{\varepsilon\to 0\\ \varepsilon>0}}<f,\varphi_{-\varepsilon}> \tag{2}$$

7.3 **Theorem** $f^-,f^+$ are distributions.

**Proof** The linearity is obvious. We give a hint on proving the fact that being given $\psi \in Diff$, the real positive sequences $(\tau_n),(\tau'_m)$ convergent to 0 and $\varphi \in \mathbf{D}$, the binary sequence:

$$< f^-, \psi \cdot \varphi_{\tau'_m} > = \lim_{n\to\infty} < f, (\psi \cdot \varphi_{\tau'_m})_{\tau_n} > \overset{not}{=} a_m \tag{3}$$

is convergent. Without loss, there exist the numbers

$$t_0 < t_1 < ... < t_{k+1} < t_{k+2} \tag{4}$$

so that

$$\varphi(t) = \varphi(t_0)\cdot\chi_{\{t_0\}}(t) \oplus \varphi(\frac{t_0+t_1}{2})\cdot\chi_{(t_0,t_1)}(t) \oplus \varphi(t_1)\cdot\chi_{\{t_1\}}(t) \oplus ...$$

$$... \oplus \varphi(\frac{t_k+t_{k+1}}{2})\cdot\chi_{(t_k,t_{k+1})}(t) \oplus \varphi(t_{k+1})\cdot\chi_{\{t_{k+1}\}}(t) \tag{5}$$

$$\psi(t) = ... \oplus \psi(t_0)\cdot\chi_{\{t_0\}}(t) \oplus \psi(\frac{t_0+t_1}{2})\cdot\chi_{(t_0,t_1)}(t) \oplus \psi(t_1)\cdot\chi_{\{t_1\}}(t) \oplus ...$$

$$... \oplus \psi(\frac{t_{k+1}+t_{k+2}}{2})\cdot\chi_{(t_{k+1},t_{k+2})}(t) \oplus \psi(t_{k+2})\cdot\chi_{\{t_{k+2}\}}(t) \oplus ... \tag{6}$$

$$\forall p \in \{0,...,k+1\}, \forall t \in (t_p, t_{p+1}), F(t_p, t_{p+1}) = F^*(t_p) \oplus F_0(t) \oplus F_*(t_{p+1}) \tag{7}$$

and the number $N \in \boldsymbol{N}$ with the property that $m,n > N$ implies:

$$t_0 < t_0 + \tau'_m < t_0 + \tau'_m + \tau_n < t_1 < ... < t_{k+1} < t_{k+1} + \tau'_m + \tau_n < t_{k+2} \tag{8}$$

The sequence $(a_m)$ is convergent, because the value of the binary numbers:

$$< f, \chi_{\{t_0+\tau'_m+\tau_n\}} >, < f, \chi_{(t_0+\tau'_m+\tau_n, t_1)} >, ...$$

$$..., < f, \chi_{\{t_{k+1}+\tau'_m\}} >, < f, \chi_{(t_{k+1}+\tau'_m, t_{k+1}+\tau'_m+\tau_n)} >, < f, \chi_{\{t_{k+1}+\tau'_m+\tau_n\}} >$$

in the hypothesis that $m,n > N$, does not depend on $m,n$ and we apply the linearity of $f$.

7.4 **Remark** The property 7.3 of $f^-, f^+$ justifies the notations 7.2 (1) and (2).

7.5 The distributions $f^-, f^+$ are called the *left limit* of $f$, respectively the *right limit* of $f$; these distributions are called together the *lateral limits* of $f$, or simply the *limits* of $f$.

7.6 The distributions that are defined by the next sums:

$$D^- f = f \oplus f^- \tag{1}$$

$$D^+ f = f \oplus f^+ \tag{2}$$

are called the *left derivative* of $f$, respectively the *right derivative* of $f$ and their common name is the *lateral derivatives,* or the *derivatives* of $f$.

7.7 All the distributions $f \in \mathbf{D}'$ admit lateral limits and all the distributions have derivatives. The next formulas of iteration of the limits and of the derivatives take place:

$$(f^-)^- = f^- \tag{1}$$

$$(f^-)^+ = f^+ \tag{2}$$

$$(f^+)^- = f^- \tag{3}$$

$$(f^+)^+ = f^+ \tag{4}$$

$$D^- D^- f = D^- f \tag{5}$$

$$D^- D^+ f = D^+ f \tag{6}$$

$$D^+D^- f = D^- f \tag{7}$$

$$D^+D^+ f = D^+ f \tag{8}$$

**Proof** We show (2):

$$< f^-,\varphi> = \lim_{\substack{\varepsilon\to 0\\ \varepsilon>0}} < f,\varphi_\varepsilon > \tag{9}$$

$$<(f^-)^+,\varphi> = \lim_{\substack{\varepsilon'\to 0\\ \varepsilon'>0}} < f^-,\varphi_{-\varepsilon'} > = \lim_{\substack{\varepsilon'\to 0\\ \varepsilon'>0}} \lim_{\substack{\varepsilon\to 0\\ \varepsilon>0}} < f,\varphi_{\varepsilon-\varepsilon'} > = \tag{10}$$

$$= \lim_{\substack{\varepsilon'\to 0\\ \varepsilon'>0}} < f,\varphi_{-\varepsilon'} > = < f^+,\varphi>$$

7.8 Taking the lateral limit, as well as the derivation are linear $\mathbf{D}'\to\mathbf{D}'$ functions (relative to the structure of $\boldsymbol{B}_2$-linear space, but not relative to the structure of $Diff$ -module of $\mathbf{D}'$).

7.9 **Examples** a) $<[\delta_{t_0}]^-,\varphi> = \lim_{\substack{\varepsilon\to 0\\ \varepsilon>0}} <[\delta_{t_0}],\varphi_\varepsilon> = \lim_{\substack{\varepsilon\to 0\\ \varepsilon>0}} \varphi_\varepsilon(t_0) =$ (1)

$$= \varphi(t_0-0) = <\delta^-_{t_0},\varphi>,\ \varphi\in\mathbf{D},\ t_0\in\boldsymbol{R}$$

thus

$$[\delta_{t_0}]^- = \delta^-_{t_0} \tag{2}$$

b) $f\in\mathbf{D}'$ is defined like at 4.10 (1). It is seen that:

$$f = f^- \tag{3}$$

7.10 Other distributions may be defined by replacing 7.1 (1), (2) with

$$\lim_{\substack{\xi\to t\\ \xi>t}} < f_{-\xi},\varphi> = \lim_{\substack{\xi\to t\\ \xi>t}} < f,\varphi_\xi > \tag{1}$$

$$\lim_{\substack{\xi\to t\\ \xi>t}} < f_{\xi},\varphi> = \lim_{\substack{\xi\to t\\ \xi>t}} < f,\varphi_{-\xi} > \tag{2}$$

etc.

## 8. The Test Functions of Two Real Variables and the Distributions that They Define

8.1 Let $\varphi_2:\boldsymbol{R}^2\to\boldsymbol{B}_2$ a function and the next conditions (that will be compared with the conditions from 2.1 a),…,e)):

a) $\varphi_2\in Diff^{(2)}$

b) the number $M>0$ exists with the property that for any $(t,u)\in\boldsymbol{R}^2$ with $\sqrt{t^2+u^2}>M$, we have:

$$\varphi_2(t,u)=0 \tag{1}$$

c) the real numbers $t',t''$ and $u',u''$ exist with $t'<t''$ and $u'<u''$ so that $supp\,\varphi_2\subset[t',t'']\times[u',u'']$, i.e. $\varphi_2$ vanishes outside a compact set

d) $\varphi_2$ is of the form:

$$\varphi_2(t,u)=\overset{k+1}{\underset{i=0}{\Xi}}\overset{p+1}{\underset{j=0}{\Xi}}\varphi_2(t_i,u_j)\cdot\chi_{\{t_i\}}(t)\cdot\chi_{\{u_j\}}(u)\oplus \tag{2}$$

$$\oplus\overset{k}{\underset{i=0}{\Xi}}\overset{p+1}{\underset{j=0}{\Xi}}\varphi_2(\frac{t_i+t_{i+1}}{2},u_j)\cdot\chi_{(t_i,t_{i+1})}(t)\cdot\chi_{\{u_j\}}(u)\oplus$$

$$\oplus\overset{k+1}{\underset{i=0}{\Xi}}\overset{p}{\underset{j=0}{\Xi}}\varphi_2(t_i,\frac{u_j+u_{j+1}}{2})\cdot\chi_{\{t_i\}}(t)\cdot\chi_{(u_j,u_{j+1})}(u)\oplus$$

$$\oplus\overset{k}{\underset{i=0}{\Xi}}\overset{p}{\underset{j=0}{\Xi}}\varphi_2(\frac{t_i+t_{i+1}}{2},\frac{u_j+u_{j+1}}{2})\cdot\chi_{(t_i,t_{i+1})}(t)\cdot\chi_{(u_j,u_{j+1})}(u)$$

where $k,p\in \boldsymbol{N}$ , $(t,u)\in \boldsymbol{R}^2$ and $t_0<t_1<...<t_{k+1}$, respectively $u_0<u_1<...<u_{p+1}$.

e) *supp* $\varphi_2$ belongs to the set ring $\boldsymbol{S}^{(2)}$ that is by definition generated by the sets:

$$(t,t')\times(u,u'), \{t\}\times(u,u'), (t,t')\times\{u\}, \{t\}\times\{u\}$$

where $t,t',u,u'\in \boldsymbol{R}$.

8.2 If $\varphi_2$ fulfills one of the next equivalent condition from 8.1:

a) and b) $\Leftrightarrow$ a) and c) $\Leftrightarrow$ d) $\Leftrightarrow$ e)

then it is called *test function* (*fundamental function*) *of two real variables*, or simply test function (fundamental function).

8.3 The space of the test functions of two real variables is also called the *fundamental space* (of the functions of two real variables) and is noted with $\mathbb{D}^{(2)}$.

8.4 From the algebraical point of view, $\mathbb{D}^{(2)}$ is a $\boldsymbol{B}_2$-algebra; $\mathbb{D}^{(2)}$ is also an ideal of $Diff^{(2)}$.

8.5 Let $\varphi_2\in\mathbb{D}^{(2)}$ a test function. We remark then that $\varphi_2(t_0,\cdot), \varphi_2(\cdot,u_0)\in\mathbb{D}$, where $t_0,u_0\in\boldsymbol{R}$.

8.6 The translation of $\varphi_2$ with $(\tau,\nu)\in\boldsymbol{R}^2$ is by definition:

$$\varphi_{2(\tau,\nu)}(t,u)=\varphi_2(t-\tau,u-\nu), t,u\in\boldsymbol{R} \tag{1}$$

8.7 We consider the quadruple $(\psi_2,(\tau_n),(\nu_m),\varphi_2)$, where $\psi_2\in Diff^{(2)},(\tau_n),(\nu_m)$ are two real positive sequences that converge to 0 and $\varphi_2\in\mathbb{D}^{(2)}$. They define the next double sequences of test functions: $(\psi_2\cdot\varphi_{2(\tau_n,\nu_m)}), (\psi_2\cdot\varphi_{2(-\tau_n,\nu_m)}), (\psi_2\cdot\varphi_{2(\tau_n,-\nu_m)})$, $(\psi_2\cdot\varphi_{2(-\tau_n,-\nu_m)})$, that we shall call *left-left convergent,..., right-right convergent.*

8.8 The $\boldsymbol{B}_2$-algebra $I_{Loc}^{(2)}$ of the *locally integrable functions* consists in the functions $f_2:\boldsymbol{R}^2\to\boldsymbol{B}_2$ with a *locally finite support*:

$$\forall t',t'',u',u''\in\boldsymbol{R}, (t',t'')\times(u',u'')\wedge supp\ f_2 \text{ is finite}$$

8.9 The *regular distributions* (over $\mathbb{D}^{(2)}$) are defined by the locally integrable functions $f_2\in I_{Loc}^{(2)}$ in the next manner:

$$< [f_2], \varphi_2 > = \pi(| supp\ f_2 \cdot \varphi_2 |) = \underset{(\xi,\xi')\in \boldsymbol{R}^2}{\Xi} f_2(\xi,\xi') \cdot \varphi_2(\xi,\xi'), \varphi_2 \in \mathbb{D}^{(2)} \quad (1)$$

8.10 We call *distribution* (over $\mathbb{D}^{(2)}$) a function $f_2 : \mathbb{D}^{(2)} \to \boldsymbol{B}_2$ that satisfies the next conditions:

a) it is linear

b) for any quadruple $(\psi_2, (\tau_n), (\nu_m), \varphi_2)$, where $\psi_2 \in Diff^{(2)}, (\tau_n), (\nu_m)$ are two real positive sequences that converge to 0 and $\varphi_2 \in \mathbb{D}^{(2)}$, the sequences:
$(f_2(\psi_2 \cdot \varphi_{2(\tau_n,\nu_m)}))$, $(f_2(\psi_2 \cdot \varphi_{2(-\tau_n,\nu_m)}))$, $(f_2(\psi_2 \cdot \varphi_{2(\tau_n,-\nu_m)}))$,
$(f_2(\psi_2 \cdot \varphi_{2(-\tau_n,-\nu_m)}))$ have a limit when $n, m \to \infty$ (for example the next property

$$\exists a \in \boldsymbol{B}_2, \exists N \in \boldsymbol{N}, \forall n,m \geq N,\ f_2(\psi_2 \cdot \varphi_{2(\tau_n,\nu_m)}) = a \quad (1)$$

is true).

8.11 The value of the distribution $f_2$ in $\varphi_2$ is noted by tradition with $< f_2, \varphi_2 >$ (instead of $f_2(\varphi_2)$).

8.12 We note with $\mathbb{D}^{(2)'}$ the $\boldsymbol{B}_2$-linear space and the $Diff^{(2)}$-module of the distributions $f_2 : \mathbb{D}^{(2)} \to \boldsymbol{B}_2$.

## 9. The Direct Product of the Distributions

9.1 a) Generally speaking, if two functions $f, g : \boldsymbol{R} \to \boldsymbol{B}_2$ are given, the following function $f \otimes g : \boldsymbol{R}^2 \to \boldsymbol{B}_2$ may be defined:

$$(f \otimes g)(t,u) = f(t) \cdot g(u), t,u \in \boldsymbol{R} \quad (1)$$

called the *direct product*, or the *tensor product* of $f$ with $g$.

b) It is true the next relation:

$$supp\ f \otimes g = supp\ f \times supp\ g \quad (2)$$

9.2 We give for 9.1 the example of two test functions $\varphi, \varphi_1 \in \mathbb{D}$; we have $\varphi \otimes \varphi_1 \in \mathbb{D}^{(2)}$, because from

$$\varphi(t) = \underset{i=0}{\overset{k+1}{\Xi}} \varphi(t_i) \cdot \chi_{\{t_i\}}(t) \oplus \underset{i=0}{\overset{k}{\Xi}} \varphi(\frac{t_i + t_{i+1}}{2}) \cdot \chi_{(t_i,t_{i+1})}(t) \quad (1)$$

$$\varphi_1(u) = \underset{j=0}{\overset{p+1}{\Xi}} \varphi_1(u_j) \cdot \chi_{\{u_j\}}(u) \oplus \underset{j=0}{\overset{p}{\Xi}} \varphi_1(\frac{u_j + u_{j+1}}{2}) \cdot \chi_{(u_j,u_{j+1})}(u) \quad (2)$$

we get:

$$(\varphi \otimes \varphi_1)(t,u) = \underset{i=0}{\overset{k+1}{\Xi}} \underset{j=0}{\overset{p+1}{\Xi}} (\varphi \otimes \varphi_1)(t_i, u_j) \cdot \chi_{\{t_i\}}(t) \cdot \chi_{\{u_j\}}(u) \oplus \quad (3)$$

$$\oplus \underset{i=0}{\overset{k}{\Xi}} \underset{j=0}{\overset{p+1}{\Xi}} (\varphi \otimes \varphi_1)(\frac{t_i + t_{i+1}}{2}, u_j) \cdot \chi_{(t_i,t_{i+1})}(t) \cdot \chi_{\{u_j\}}(u) \oplus$$

$$\oplus \underset{i=0}{\overset{k+1}{\Xi}} \underset{j=0}{\overset{p}{\Xi}} (\varphi \otimes \varphi_1)(t_i, \frac{u_j + u_{j+1}}{2}) \cdot \chi_{\{t_i\}}(t) \cdot \chi_{(u_j,u_{j+1})}(u) \oplus$$

$$\oplus \underset{i=0}{\overset{k}{\Xi}} \underset{j=0}{\overset{p}{\Xi}} (\varphi \otimes \varphi_1)(\frac{t_i + t_{i+1}}{2}, \frac{u_j + u_{j+1}}{2}) \cdot \chi_{(t_i, t_{i+1})}(t) \cdot \chi_{(u_j, u_{j+1})}(u)$$

In (1),...,(3), we have supposed that $k, p \in \boldsymbol{N}, t, u \in \boldsymbol{R}, t_0 < t_1 < ... < t_{k+1}$ and $u_0 < u_1 < ... < u_{p+1}$.

9.3 Our purpose in this paragraph is that of transferring the definition of the direct product of the functions to distributions. For $f, g \in I_{Loc}$ we have $f \otimes g \in I_{Loc}^{(2)}, [f \otimes g] \in \mathbb{D}^{(2)'}$ and, by using the same notation ' $\otimes$ ' for the direct product of the distributions to be defined, we must have:

$$[f \otimes g] = [f] \otimes [g] \qquad (1)$$

9.4 **Theorem** Let $f \in \mathbb{D}'$ and $\varphi_2 \in \mathbb{D}^{(2)}$. The function $\varphi' : \boldsymbol{R} \to \boldsymbol{B}_2$ that is defined in the following manner:

$$\varphi'(t) = < f, \varphi_2(t, \cdot) >, t \in \boldsymbol{R} \qquad (1)$$

is a test function.

**Proof** a) Let us remark in the beginning that because for any fixed $t \in \boldsymbol{R}$, the function $\varphi_2(t, \cdot)$ is a test function (see 8.5), the formula $< f, \varphi_2(t, \cdot) >$ that defines $\varphi'(t)$ makes sense.

b) We refer to the expression 8.1 (2) of $\varphi_2$ and we see by a direct computation that:

$$\varphi'(t) = \underset{i=0}{\overset{k+1}{\Xi}} \varphi'(t_i) \cdot \chi_{\{t_i\}}(t) \oplus \underset{i=0}{\overset{k}{\Xi}} \varphi'(\frac{t_i + t_{i+1}}{2}) \cdot \chi_{(t_i, t_{i+1})}(t) \qquad (2)$$

where

$$\varphi'(t_i) = \underset{j=0}{\overset{p+1}{\Xi}} \varphi_2(t_i, u_j) \cdot F_0(u_j) \oplus \underset{j=0}{\overset{p}{\Xi}} \varphi_2(t_i, \frac{u_j + u_{j+1}}{2}) \cdot F(u_j, u_{j+1}), i = \overline{0, k+1} \qquad (3)$$

$$\varphi'(\frac{t_i + t_{i+1}}{2}) = \underset{j=0}{\overset{p+1}{\Xi}} \varphi_2(\frac{t_i + t_{i+1}}{2}, u_j) \cdot F_0(u_j) \oplus$$

$$\oplus \underset{j=0}{\overset{p}{\Xi}} \varphi_2(\frac{t_i + t_{i+1}}{2}, \frac{u_j + u_{j+1}}{2}) \cdot F(u_j, u_{j+1}), i = \overline{0, k} \qquad (4)$$

$\varphi'$ is a test function.

9.5 **Notation** Another manner of writing the relation 9.4 (1) is:

$$\varphi'(t) =< f(u), \varphi_2(t, u) >, t, u \in \boldsymbol{R} \qquad (1)$$

9.6 **Remark** The advantage of the notation 9.5 is that of indicating that $f$ acts on the second variable of $\varphi_2$. The disadvantage is its abusive character, the argument of $f$ is $\varphi_2(t, \cdot) \in \mathbb{D}$, not $u \in \boldsymbol{R}$.

When we work with functions that have several variables, notations like 9.5 can be useful.

9.7 **Theorem** It is true:

$$D_t^- \varphi'(t) = < f, D_t^- \varphi_2(t, \cdot) >, t \in \boldsymbol{R} \qquad (1)$$

$$D_t^+ \varphi'(t) = < f, D_t^+ \varphi_2(t, \cdot) >, t \in \boldsymbol{R} \qquad (2)$$

where the symbols $D_t^-, D_t^+$ refer to the lateral derivation relative to $t$.

**Proof** The relation

$$\varphi'(t+\tau) =< f, \varphi_2(t+\tau,\cdot) >, t \in \boldsymbol{R} \quad (3)$$

may be written for any $\tau \in \boldsymbol{R}$ and, as $\varphi'$ is a test function, we may take the limit in (3), for example:

$$\varphi'(t-0) =< f, \varphi_2(t-0,\cdot) >, t \in \boldsymbol{R} \quad (4)$$

By the summation of (4) with 9.4 (1), 9.7 (1) results.

9.8 Let $f, g \in \mathbf{D}', \varphi_2 \in \mathbf{D}^{(2)}$ and $\varphi' \in \mathbf{D}$ defined by 9.4 (1). When $\varphi_2$ is variable (and $\varphi'$ is also variable) $f$ and $g$ being fixed, the formula:

$$< f \otimes g, \varphi_2 > = < g, \varphi' > \quad (1)$$

can also be abusively written under the form

$$< (f \otimes g)(t,u), \varphi_2(t,u) > = < g(t), \varphi'(t) > = < g(t), < f(u), \varphi_2(t,u) >>, t,u \in \boldsymbol{R} \quad (2)$$

It defines a function $f \otimes g : \mathbf{D}^{(2)} \to \boldsymbol{B}_2$ that is called the *direct product* (or the *tensor product*) of the distributions $f$ and $g$.

9.9 From the formulas 9.4 (2), (3), (4) and with the notations $G, G_0$ for the fundamental functions that are associated to $g$, we have the next computation of the value of $f \otimes g$ in $\varphi_2$:

$$< f \otimes g, \varphi_2 > = \underset{i=0}{\overset{k+1}{\Xi}} \underset{j=0}{\overset{p+1}{\Xi}} \varphi_2(t_i, u_j) \cdot G_0(t_i) \cdot F_0(u_j) \oplus \quad (1)$$

$$\oplus \underset{i=0}{\overset{k}{\Xi}} \underset{j=0}{\overset{p+1}{\Xi}} \varphi_2(\frac{t_i + t_{i+1}}{2}, u_j) \cdot G(t_i, t_{i+1}) \cdot F_0(u_j) \oplus$$

$$\oplus \underset{i=0}{\overset{k+1}{\Xi}} \underset{j=0}{\overset{p}{\Xi}} \varphi_2(t_i, \frac{u_j + u_{j+1}}{2}) \cdot G_0(t_i) \cdot F(u_j, u_{j+1}) \oplus$$

$$\oplus \underset{i=0}{\overset{k}{\Xi}} \underset{j=0}{\overset{p}{\Xi}} \varphi_2(\frac{t_i + t_{i+1}}{2}, \frac{u_j + u_{j+1}}{2}) \cdot G(t_i, t_{i+1}) \cdot F(u_j, u_{j+1})$$

9.10 **Theorem** $f \otimes g \in \mathbf{D}^{(2)}{}'$.

**Proof** The linearity is obvious and the property of convergence 8.10 b) results from the linearity, from the formula 9.9 (1) and from the properties of the fundamental functions $F$, $F_0, G, G_0$.

9.11 **Theorem** Let $f, g, h \in \mathbf{D}'$. The next properties are true:

a) $$f \otimes g = g \otimes f \quad (1)$$

in the sense that the next relation is fulfilled:

$$\forall \varphi_2 \in \mathbf{D}^{(2)}, < g(t), < f(u), \varphi_2(t,u) >> = < f(u), < g(t), \varphi_2(t,u) >> \quad (2)$$

b) $$(f \oplus g) \otimes h = (f \otimes h) \oplus (g \otimes h) \quad (3)$$

c) After defining the spaces of test functions $\mathbf{D}^{(3)}$ and of the distributions $\mathbf{D}^{(3)} \to \boldsymbol{B}_2$, we have the associativity:

$$(f \otimes g) \otimes h = f \otimes (g \otimes h) \quad (4)$$

**Proof** a) It results from 9.9.

b) Obvious.

c) It is necessary to replace the formula 9.9 (1) with the formula:

$$< f \otimes g \otimes h, \varphi_3 > = \Xi_{i=0}^{k+1} \Xi_{j=0}^{p+1} \Xi_{s=0}^{l+1} \varphi_3(t_i, u_j, v_s) \cdot H_0(t_i) \cdot G_0(u_j) \cdot F_0(v_s) \oplus ... \quad (5)$$

$$... \oplus \Xi_{i=0}^{k} \Xi_{j=0}^{p} \Xi_{s=0}^{l} \varphi_3\left(\frac{t_i + t_{i+1}}{2}, \frac{u_j + u_{j+1}}{2}, \frac{v_s + v_{s+1}}{2}\right) \cdot H(t_i, t_{i+1}) \cdot G(u_j, u_{j+1}) \cdot F(v_s, v_{s+1})$$

where $\varphi_3 \in \mathbb{D}^{(3)}$ etc.

9.12 Let $\psi \in Diff$, $f, g \in \mathbb{D}'$, $\varphi_2 \in \mathbb{D}^{(2)}$. The fact that:

$$< \psi(u) \cdot f(u), < g(t), \varphi_2(t,u) >> = < f(u), \psi(u) \cdot < g(t), \varphi_2(t,u) >> = \quad (1)$$

$$= < f(u), < g(t), \psi(u) \cdot \varphi_2(t,u) >> = < \psi \cdot (f \otimes g), \varphi_2 >$$

where $t, u \in \boldsymbol{R}$ allows us to write briefly

$$(\psi \cdot f) \otimes g = \psi \cdot (f \otimes g) \quad (2)$$

Anyway, the relation (2) must be understood keeping in mind the fact that the 'arguments' of $\psi, f, g$ are $u, u, t$.

9.13 For $f, g$ like before, we study the behavior of $f \otimes g$ relative to the derivation of the distributions. Because the function

$$\psi(u) =< g, \varphi_2(\cdot, u) >, u \in \boldsymbol{R} \quad (1)$$

is a test function we can write

$$< D_u^-(f \otimes g)(t,u), \varphi_2(t,u) > = \quad (2)$$

$$= < (f \otimes g)(t,u), \varphi_2(t,u) > \oplus \lim_{\substack{\varepsilon \to 0 \\ \varepsilon > 0}} < (f \otimes g)(t,u), \varphi_2(t, u - \varepsilon) > =$$

$$= < f(u), < g(t), \varphi_2(t,u) >> \oplus \lim_{\substack{\varepsilon \to 0 \\ \varepsilon > 0}} < f(u), < g(t), \varphi_2(t, u - \varepsilon) >> =$$

$$= < f(u), < g(t), \varphi_2(t,u) >> \oplus \lim_{\substack{\varepsilon \to 0 \\ \varepsilon > 0}} < f(u + \varepsilon), < g(t), \varphi_2(t,u) >> =$$

$$= \lim_{\substack{\varepsilon \to 0 \\ \varepsilon > 0}} < f(u) \oplus f(u + \varepsilon), < g(t), \varphi_2(t,u) >> = < ((D_u^- f) \otimes g)(t,u), \varphi_2(t,u) >$$

If we keep in mind the variables $u, u, t$ to which the left derivation operator, $f$ and $g$ refer, the relation (2) can be briefly written like this:

$$D^-(f \otimes g) = (D^- f) \otimes g \quad (3)$$

9.14 **Examples** of direct product of distributions. a) Let $\{t_z \mid z \in \boldsymbol{Z}\}, \{u_w \mid w \in \boldsymbol{Z}\}$ two locally finite families and $[\Xi_{z \in \boldsymbol{Z}} \delta_{t_z}], [\Xi_{w \in \boldsymbol{Z}} \delta_{u_w}] \in \mathbb{D}'$ the regular distributions that they define. We have for $\varphi_2 \in \mathbb{D}^{(2)}$:

$$<[\Xi_{w \in \boldsymbol{Z}} \delta_{u_w}] \otimes [\Xi_{z \in \boldsymbol{Z}} \delta_{t_z}], \varphi_2 > = \int_{-\infty}^{\infty} {}_t \Xi_{z \in \boldsymbol{Z}} \delta(t - t_z) \int_{-\infty}^{\infty} {}_u \Xi_{w \in \boldsymbol{Z}} \delta(u - u_w) \cdot \varphi_2(t,u) = \quad (1)$$

$$= \Xi_{z \in \boldsymbol{Z}} \Xi_{w \in \boldsymbol{Z}} \varphi_2(t_z, u_w)$$

where the symbols $\int_{-\infty}^{\infty} {}_t, \int_{-\infty}^{\infty} {}_u$ indicate the variables relative to which the integration is done.

b) We note with $g \in \mathbf{D}'$ the distribution 5.7 (1) and we have for $\varphi_2 \in \mathbf{D}^{(2)}$ and $u_0 \in \boldsymbol{R}$:

$$< \delta_{u_0}^{-} \otimes g, \varphi_2 > = < g(t), < \delta_{u_0}^{-}, \varphi_2(t,u) >> = \int_{-\infty}^{\infty} D^{-}\varphi_2(\cdot, u_0 - 0) \quad (2)$$

## 10. Distributions Over Test Functions of Several Variables

10.1 This generalization refers to the number of independent variables of the test functions, in the sense of the definition of the spaces $\mathbf{D}, \mathbf{D}^{(2)}, ..., \mathbf{D}', \mathbf{D}^{(2)}{}', ...$, problem that was already dealt with in the previous paragraphs. We give a hint on the steps to be followed.
10.2 Step I, the definition of the differentiability of the functions $\varphi : \boldsymbol{R}^n \to \boldsymbol{B}_2$ (by generalizing 1.12, 1.17).
10.3 Step II, we define the test functions $\varphi \in \mathbf{D}^{(n)}$ to be the differentiable functions $\varphi \in Diff^{(n)}$ with the property that they have a bounded support.
10.4 Step III, we define the distributions $f \in \mathbf{D}^{(n)}{}'$ to be the linear functions $f : \mathbf{D}^{(n)} \to \boldsymbol{B}_2$ with the property that for any $n+2$-tuple $(\psi, (\tau_{p_1}), ..., (\tau_{p_n}), \varphi)$, where $\varphi \in Diff^{(n)}$, $(\tau_{p_1}), \ldots, (\tau_{p_n})$ are real positive sequences that converge to zero and $\varphi \in \mathbf{D}^{(n)}$ we have that $\lim\limits_{p_1, ..., p_n \to \infty} f(\psi \cdot \varphi_{(\pm\tau_{p_1}, ..., \pm\tau_{p_n})})$ exists (for example the next property

$$\exists a \in \boldsymbol{B}_2, \exists N \in \boldsymbol{N}, \forall p_1, ..., p_n \geq N, f(\psi \cdot \varphi_{(\tau_{p_1}, ..., \tau_{p_n})}) = a \quad (1)$$

is true).

We have noted

$$\varphi_{(\pm\tau_{p_1}, ..., \pm\tau_{p_n})}(t_1, ..., t_n) = \varphi(t_1 \mp \tau_{p_1}, ..., t_n \mp \tau_{p_n}),\ t_1, ..., t_n \in \boldsymbol{R} \quad (2)$$

## 11. Other Distributions

11.1 As opposed to 10, this second type of generalization of the paragraphs 2,…,9 refers to the functions $\varphi : \boldsymbol{R} \to \boldsymbol{B}_2$ and starts from redefining the notion of test function.

11.2 a) It is called *space of test functions* (or *fundamental space*) a set $\mathbf{K} \subset Diff$ with the next properties:

a.1) $\varphi \in \mathbf{K}, \tau \in \boldsymbol{R} \Rightarrow \varphi_\tau \in \mathbf{K}$

a.2) $\varphi \in \mathbf{K} \Rightarrow \varphi^{-}, \varphi^{+} \in \mathbf{K}$

a.3) $\varphi, \varphi' \in \mathbf{K} \Rightarrow \varphi \oplus \varphi' \in \mathbf{K}$

a.4) $\psi \in Diff, \varphi \in \mathbf{K} \Rightarrow \psi \cdot \varphi \in \mathbf{K}$

b) The functions $\varphi \in \mathbf{K}$ are called *test functions* (or *fundamental functions*).

11.3 **Remarks** a) All these properties were used in the previous paragraphs.

b) $\mathbf{K}$ is invariant to translations, it contains a function together with its lateral limits, it is a $\boldsymbol{B}_2$-algebra and it is also an ideal of $Diff$.

11.4 **Examples** of spaces of test functions. $\{0\}$ and $Diff$ are the two improper subspaces of $Diff$ that fulfill the conditions of the definition 11.2. We have noted with $\{0\}$ the set formed by the null function $\boldsymbol{R} \to \boldsymbol{B}_2$.

Other examples are $\mathbb{D}$, $I_{Loc}$ and $I_\infty$.

11.5 **Proposition** If $\mathbb{K} \neq \{0\}$ is a space of test functions, then $\chi_{\{t_0\}} \in \mathbb{K}, t_0 \in \boldsymbol{R}$. $I_\infty$ is the smallest non-null space of test functions (in the sense of the order given by the inclusion).

**Proof** Let $\varphi \in \mathbb{K}$ a non-null test function (we can take such a function because $\mathbb{K} \neq \{0\}$) and let us suppose that $t_0$ is a point where $\varphi$ is equal with 1. Because $\chi_{\{t_0\}} \in Diff$, from 11.2 a.4) we have that

$$\chi_{\{t_0\}} \cdot \varphi = \chi_{\{t_0\}} \in \mathbb{K} \tag{1}$$

11.6 Let $\mathbb{K}$ a space of test functions. We call *distribution* (over $K$) a linear function $f : \mathbb{K} \to \boldsymbol{B}_2$ that fulfills one of the equivalent conditions:

a) for any triple $(\psi, (\tau_n), \varphi)$ with $\psi \in Diff$, $(\tau_n)$ a real positive sequence that converges to 0 and $\varphi \in \mathbb{K}$, the binary sequences $(f(\psi \cdot \varphi_{\tau_n})), (f(\psi \cdot \varphi_{-\tau_n}))$ are convergent.

b) for any couple $(\psi, \varphi)$ with $\psi \in Diff, \varphi \in \mathbb{K}$, the function in $t$: $f(\psi \cdot \varphi_t)$ is differentiable.

11.7 a) The distribution $f$ over $\mathbb{K}$ is said to be *regular*, or *of function type*, if a function $g : \boldsymbol{R} \to \boldsymbol{B}_2$ exists with the property that

$$f(\varphi) = \int_{-\infty}^{\infty} g \cdot \varphi, \varphi \in \mathbb{K} \tag{1}$$

b) If $f$ is not regular, then it is called *singular*.

11.8 By tradition, we shall note with $< f, \varphi >$ (instead of $f(\varphi)$) the value of the distribution $f$ in the 'point' $\varphi$; other notations are $< f(t), \varphi(t) >, (f, \varphi)$ etc.

If the distribution is regular, we shall continue to use the notation $< [g], \varphi >$.

11.9 The set $\mathbb{K}'$ of the $\mathbb{K} \to \boldsymbol{B}_2$ distributions is obviously organized as a $\boldsymbol{B}_2$-linear space and $Diff$-module. The translation of a distribution, the limit and the derivative are easily defined.

11.10 We have the next inclusions of spaces of test functions and of spaces of distributions:

$$\{0\} \subset I_\infty \subset I_{Loc} \subset Diff, \{0\}' \supset I_\infty' \supset I_{Loc}' \supset Diff'$$

$$I_\infty \subset \mathbb{D} \subset Diff, I_\infty' \supset \mathbb{D}' \supset Diff'$$

$$\{0\} \subset \mathbb{K} \subset Diff, \{0\}' \supset \mathbb{K}' \supset Diff'$$

where $\mathbb{K}$ is an arbitrary space of test functions.

11.11 **Notation** For $\mathbb{K}$ like before we note:

$$Reg\,\mathbb{K} = \{f \mid f : \boldsymbol{R} \to \boldsymbol{B}_2, [f] \in \mathbb{K}'\} \tag{1}$$

11.12 $Reg\,\mathbb{K}$ is a $\boldsymbol{B}_2$-linear subspace and a $Diff$-submodule of $Diff$.

11.13 **Theorem** The next statements are true:

a) $Reg\,\{0\} = \boldsymbol{B}_2^{\boldsymbol{R}}$ a)' $\{0\}' = \{[f] \mid f \in \boldsymbol{B}_2^{\boldsymbol{R}}\} = \{0 \mid 0 : \{0\} \to \boldsymbol{B}_2\}$ (1)

b.1) $Reg\ I_\infty = Diff$ b.1)' $I_\infty^{'} = \{[f] \mid f \in Diff\}$ (2)

b.2) $Reg\ Diff = I_\infty$ (3)

c.1) $Reg\ I_{Loc} = \mathbf{D}$ (4)

c.2) $Reg\ \mathbf{D} = I_{Loc}$ (5)

**Proof** a) Let $f \in \boldsymbol{B}_2^{\boldsymbol{R}}$ arbitrary. It defines the null function:

$$\forall \varphi \in \{0\}, \int_{-\infty}^{\infty} f \cdot \varphi = 0 \tag{6}$$

that is a distribution of $\{0\}'$.

a)' The only linear function $\{0\} \to \boldsymbol{B}_2$ is the null function, that is a distribution.

b.1) Let $f : \boldsymbol{R} \to \boldsymbol{B}_2$ arbitrary. From the fact that $\chi_{\{t\}} \in I_\infty$, where $t \in \boldsymbol{R}$, the condition

$$< [f], \chi_{\{t\}} > = f(t) \in Diff \tag{7}$$

shows that $Reg\ I_\infty \subset Diff$. The inverse inclusion is obvious.

b.1)' Let $f \in I_\infty^{'}$, for which we have:

$$F_0(t) \overset{def}{=} < f, \chi_{\{t\}} > \in Diff \tag{8}$$

Any non-null test function from $I_\infty$ is of the form $\chi_{\{t_0\}} \oplus ... \oplus \chi_{\{t_n\}}$ and $f$ acts on it in the next manner:

$$< f, \chi_{\{t_0\}} \oplus ... \oplus \chi_{\{t_n\}} > = < f, \chi_{\{t_0\}} > \oplus ... \oplus < f, \chi_{\{t_n\}} > = F_0(t_0) \oplus ... \oplus F_0(t_n) = \tag{9}$$

$$= < [F_0], \chi_{\{t_0\}} \oplus ... \oplus \chi_{\{t_n\}} >$$

i.e. $f$ is of function type:

$$f = [F_0] \tag{10}$$

showing that $I_\infty^{'} \subset \{[f] \mid f \in Diff\}$. The inverse inclusion results from b.1).

b.2) Let $f : \boldsymbol{R} \to \boldsymbol{B}_2$. Because $1 \in Diff$ and $supp\ f \cdot 1$ is finite, we infer that $f \in I_\infty$, thus $Reg\ Diff \subset I_\infty$. The fact that $I_\infty \subset Reg\ Diff$ is easily proved.

c.1) Let the function $f : \boldsymbol{R} \to \boldsymbol{B}_2$. As $\chi_{\{t\}} \in I_{Loc}$, like at (7) $f$ defines a regular distribution if it is differentiable and moreover if $\forall \varphi \in I_{Loc}, supp\ f \cdot \varphi$ is finite i.e. if $supp\ f$ is bounded, resulting that $f \in \mathbf{D}$. We have proved that $Reg\ I_{Loc} \subset \mathbf{D}$. The inverse inclusion is obvious.

c.2) Let $f : \boldsymbol{R} \to \boldsymbol{B}_2$. The condition $\forall \varphi \in \mathbf{D}, supp\ f \cdot \varphi$ is finite implies $f \in I_{Loc}$, thus $Reg\ \mathbf{D} \subset I_{Loc}$. The inverse inclusion is obvious, because this was the definition of the regular distributions.

11.14 We have here the open problem: does $I_{Loc}^{'}$ contain singular distributions ?

## 12. The Convolution Algebras of the Distributions From $\mathbf{D}'$

12.1 Let $f, g \in \mathbf{D}'$ two distributions. The function

$$\psi(u) = < g, \varphi_{-u} >, u \in \boldsymbol{R}, \varphi \in \mathbf{D} \tag{1}$$

is a differentiable function (with an unbounded support, generally), because it is a sum of differentiable functions of the form $G_0(t_i - u), i = \overline{0, k+1}, G(t_i - u, t_{i+1} - u), i = \overline{0, k}, G_0, G$ being the fundamental functions of $g$ and this is why it makes sense sometimes, for example if $f \in Diff' \subset \mathbf{D}'$, to refer to $< f, \psi >$.

12.2 It is called the *convolution product* of the distributions $f$ and $g$ (in this order) the function $f * g : \mathbf{D} \to \boldsymbol{B}_2$ that is defined by

$$< f * g, \varphi > = < f, \psi > = < f(u), < g(t), \varphi(t+u) >> \quad (1)$$

where $\varphi \in \mathbf{D}$ and $t, u \in \boldsymbol{R}$.

12.3 As resulted by the successive application of two distributions, the convolution product is a distribution, this is obvious. This fact justifies the notation that was used at 12.2.

12.4 It is a natural desire to present the convolution product of the distributions to be their direct product, under the form:

$$< f * g, \varphi > = < f \otimes g, \varphi \circ \sigma > \quad (1)$$

where $\sigma : \boldsymbol{R}^2 \to \boldsymbol{R}$ is the sum function:

$$\sigma(t, u) = t + u, t, u \in \boldsymbol{R} \quad (2)$$

This fact is not possible however, because when $\varphi \in \mathbf{D}, \varphi \circ \sigma \notin \mathbf{D}^{(2)}$; $\varphi \circ \sigma$ is not, generally, differentiable or with bounded support.

12.5 **Counterexample** at 12.4. The function $\delta \in \mathbf{D}$ fulfills the previous property: the function of two variables

$$(\delta \circ \sigma)(t, u) = \delta(t + u), t, u \in \boldsymbol{R} \quad (1)$$

does not belong to $Diff^{(2)}$ and it has also an unbounded support.

12.6 Let $f, g : \boldsymbol{R} \to \boldsymbol{B}_2$; their convolution product is the function $f * g : \boldsymbol{R} \to \boldsymbol{B}_2$ defined by

$$(f * g)(t) = \int_{-\infty}^{\infty} \xi\, f(\xi) \cdot g(t - \xi) \quad (1)$$

12.7 We define the following sets of locally integrable functions $f : \boldsymbol{R} \to \boldsymbol{B}_2$:

$$I_{Inf} = \{ f \mid \forall \alpha \in \boldsymbol{R}, (-\infty, \alpha) \wedge supp\ f\ is\ finite \} \quad (1)$$

$$I_{Sup} = \{ f \mid \forall \alpha \in \boldsymbol{R}, (\alpha, \infty) \wedge supp\ f\ is\ finite \} \quad (2)$$

12.8 In the previous definition, a set $H \subset \boldsymbol{R}$ satisfying

$$\forall \alpha \in \boldsymbol{R}, (-\infty, \alpha) \wedge H \text{ is finite}$$

$$\forall \alpha \in \boldsymbol{R}, (\alpha, \infty) \wedge H \text{ is finite}$$

is called *inferiorly finite*, respectively *superiorly finite*.

12.9 There exists the situation when the convolution product of $f, g \in I_{Loc}$ has sense, for example if $f, g$ belong to $I_\infty, I_{Inf}, I_{Sup}$ and in all these cases $f * g \in I_{Loc}$ (more exactly: $f * g$ belongs to $I_\infty, I_{Inf}$, respectively $I_{Sup}$). In such situations $[f * g] \in D'$ and we can write that:

$$<[f * g], \varphi > = \underset{\xi \in \boldsymbol{R}}{\Xi} (f * g)(\xi)\varphi(\xi) = \quad (1)$$

$$= \underset{\xi \in \boldsymbol{R}}{\Xi} \underset{u \in \boldsymbol{R}}{\Xi} f(u) \cdot g(\xi - u) \cdot \varphi(\xi) = \underset{u \in \boldsymbol{R}}{\Xi} \underset{\xi \in \boldsymbol{R}}{\Xi} f(u) \cdot g(\xi - u) \cdot \varphi(\xi)$$

With the substitution

$$\xi - u = t \qquad (2)$$

that transforms the summation variables $\xi, u$ in the summation variables $t, u$, it results

$$<[f * g], \varphi> = \underset{u \in \mathbf{R}}{\Xi} \; \underset{t \in \mathbf{R}}{\Xi} f(u) \cdot g(t) \cdot \varphi(t+u) = <[f](u), <[g](t), \varphi(t+u)>> \qquad (3)$$

12.10 The formula 12.9 (3) rediscovers the definition 12.2 (1) for the case of the regular distributions. By comparing these formulas we have that

$$[f * g] = [f] * [g] \qquad (1)$$

12.11 The regular distribution $[\delta]$ has a convolution product with any $f \in \mathbf{D}'$ and we can write that

$$<[\delta] * f, \varphi> = <[\delta](u), < f(t), \varphi(t+u)>> = < f(t), \varphi(t)> = \qquad (1)$$

$$= < f(t), <[\delta](u), \varphi(t+u)>> = < f * [\delta], \varphi>, \varphi \in \mathbf{D}$$

12.12 It is said that a linear subspace $\mathbf{U}$ of $\mathbf{D}'$ is a *convolution algebra* (of distributions) if it has the next properties:

a) $[\delta] \in \mathbf{U}$

b) the convolution product of two or several distributions from $\mathbf{U}$ is defined and it belongs to $\mathbf{U}$. The convolution product is associative and commutative.

12.13 **Examples** of convolution algebras.

a) $\{0, [\delta]\}$ is the simplest convolution algebra, where 0 is the null distribution.

b) $\{0, [\delta], \delta^-, \delta^+, [\delta] \oplus \delta^-, [\delta] \oplus \delta^+, \delta^- \oplus \delta^+, [\delta] \oplus \delta^- \oplus \delta^+\}$ is a convolution algebra that has the remarkable property that it contains a distribution $f$ together with the distributions $f^-$ and $f^+$. This fact is also true at the next example.

c) Let $I, J, K \subset \mathbf{R}$ finite sets (inferiorly finite, superiorly finite sets). The space of the distributions of the form

$$\underset{i \in I}{\Xi} [\delta_i] \oplus \underset{j \in J}{\Xi} [\delta_j^-] \oplus \underset{k \in K}{\Xi} [\delta_k^+]$$

is a convolution algebra.

12.14 **Proposition** Let the convolution algebra $\mathbf{U}$. The next statements are true:

a) $$f * (g \oplus h) = (f * g) \oplus (f * h) \qquad (1)$$

b) If $f^-, f^+ \in \mathbf{U}$, then $D^- f * g, f * D^- g, D^-(f * g) \in \mathbf{U}$ and the next result holds:

$$D^-(f * g) = D^- f * g = f * D^- g \qquad (2)$$

At a), b) the distributions $f, g, h$ belong to $\mathbf{U}$.

**Proof** b) The fact that $D^- f, D^- g, D^- f * g, f * D^- g, D^-(f * g) \in \mathbf{U}$ is obvious. We have:

$$< D^-(f * g), \varphi> = < f * g, \varphi> \oplus \lim_{\substack{\varepsilon \to 0 \\ \varepsilon > 0}} < f * g, \varphi_\varepsilon> = \qquad (3)$$

$$= < f, \psi> \oplus \lim_{\substack{\varepsilon \to 0 \\ \varepsilon > 0}} < f, \psi_\varepsilon> = < D^- f, \psi> = < D^- f * g, \varphi>$$

where $\varphi \in \mathbf{D}$ and $\psi \in Diff$ is given by 12.1 (1).